\documentclass[reqno,a4paper,11pt,oneside]{amsart}

\usepackage{algpseudocode,algorithm,algorithmicx}
\usepackage{mathtools,amssymb}
\usepackage{todonotes}
\usepackage{booktabs}
\usepackage{multirow}
\usepackage{graphicx,epstopdf} 
\usepackage[caption=false]{subfig}
\usepackage[final,stretch=10]{microtype}
\usepackage[hidelinks,pdftex]{hyperref}

\newcommand{\B}{T} %supremizer operator
\newcommand{\E}{E}
\newcommand{\p}{\mathcal{P}}
\newcommand{\I}{\mathbb{I}}

\newcommand{\parmu}{\boldsymbol\mu}
\newcommand{\vol}{\nu}

\newcommand{\NB}{{\rm N}}%boundary Dir
\newcommand{\DB}{{\rm D}}%boundary Neum
\newcommand{\pf}{\chi} %log transformed pay of
\newcommand{\g}{g}% payoff in sp formulation
\newcommand{\dual}{\chi}% dual basis functions
\newcommand{\Nopt}{M}% number of options
\newcommand{\popt}{\p_{\rm opt}}% number of options
\newcommand{\tdeam}{D}% number of options
\newcommand{\muopti}{\Theta}
\newcommand{\muopt}{\boldsymbol{\muopti}}
\newcommand{\muoptstar}{\muopt^{\star}}
\newcommand{\LL}{{L^2(\Omega)}} %l2 product
\newcommand{\V}{{V^{\prime}}}

\DeclareMathOperator{\spann}{span}
\DeclareMathOperator*{\argmin}{arg\,min}

\DeclareMathOperator{\train}{train}

\DeclareMathOperator{\tol}{tol}

\newcommand{\N}{\mathcal{N}}

\usepackage{ushort}
\newcommand*{\Mat}[1]{{\mathbf{#1}}}

\newcommand*{\deam}[1]{{\widetilde{#1}}}
\numberwithin{equation}{section}

\newtheorem{remark}{Remark}[section]
\begin{document}

\title{Model reduction for calibration of American options}

\author{O. Burkovska
\and K. Glau
\and M. Mahlstedt
\and B. Wohlmuth}

\address{O. Burkovska, Institute for Numerical Mathematics, Technische Universit\"at M\"unchen, 85748
  Garching b.~M\"unchen, Germany}
\email{burkovsk@ma.tum.de}
\address{K. Glau \\
  Chair of Mathematical Finance \\
  Technische Universit\"at M\"unchen \\
  85748 Garching b.~M\"unchen, Germany}
\email{kathrin.glau@tum.de}
\address{M. Mahlstedt \\
  Chair of Mathematical Finance \\
  Technische Universit\"at M\"unchen \\
  85748 Garching b.~M\"unchen, Germany}
\email{mirco.mahlstedt@tum.de}
\address{B. Wohlmuth \\
  Institute for Numerical Mathematics \\
  Technische Universit\"at M\"unchen \\
  85748 Garching b.~M\"unchen, Germany}
\email{wohlmuth@ma.tum.de}

\thanks{This work was partly supported by: DFG grant WO671/11-1; International Research
  Training Group IGDK1754, funded by the German Research Foundation (DFG) and the Austrian
  research fund (FWF); KPMG Center of Excellence in Risk Management}

\date{\today}


\begin{abstract}
 American put options are among the most frequently traded single stock options, and their
calibration is computationally challenging since no closed-form expression is available.
Due to the higher flexibility in comparison to European options, the mathematical model involves additional constraints, and a variational inequality is obtained.
We use the Heston stochastic volatility model to describe the price of a single stock option.
In order to speed up the calibration process, we apply  two  model
reduction strategies. Firstly,  a reduced basis method (RBM) is used to define a suitable low-dimensional basis
for the numerical approximation of the parameter-dependent partial differential equation ($\mu$PDE)  model. By doing so the computational complexity for solving
the $\mu$PDE is drastically reduced, and applications of standard minimization algorithms for the calibration are significantly faster than working with a high-dimensional 
finite element basis. Secondly,  so-called de-Americanization strategies are applied. Here, the main idea is to reformulate the  calibration problem for American options as a problem for European options and to exploit
closed-form solutions.  Both reduction techniques are systematically compared and tested for both synthetic and market data sets.
\end{abstract}

\maketitle

\section{Introduction}\label{ch6:sec:motivation}
Mathematical models for option pricing typically depend sensitively on a set of parameters such as, e.g., the interest rate, the long-run variance, the  rate of mean reversion, the volatility of volatility and the correlation parameter. Reliable predictions are only possible if these parameters are known. While the interest rate is often known a priori, other parameters cannot be accessed directly from the observed market data but have to be determined by a computationally challenging calibration process. Constantly changing market situations then require fast parameter fitting algorithms.
Here, we consider
single stock options of American type. This type of  options is path-dependent, as the option holder has the right to exercise it at any time until maturity. Moreover, it is one of the most popular choices to be traded at exchange stocks.
To be able to react instantaneously to market movements, complexity reduction techniques are of special interest. Both European and American options can be based on a parameter-dependent partial differential equation in time and asset price, {see, e.g.,~\cite{achdou,hilber} and the references therein}.
However due to the higher flexibility of exercising American options compared to European options, the mathematical model for an American model has to be enriched by a suitable inequality constraint, reflecting the arbitrage-free principle.
Then the problem under consideration can be reformulated as a weak variational inequality problem to which semi-smooth Newton schemes can be applied as non-linear solvers. Classical finite discretization schemes such as finite elements or finite differences then require rather high-dimensional basis spaces leading to large systems to be solved. To reduce the computational cost, we 
apply two different  reduction techniques. The first is based on reduced basis methods, which are perfectly suitable for 
parameter-dependent partial differential equations, while the second exploits the fact that closed-form solutions can be accessed for the simpler European options.
The key idea  of the reduced basis method is to replace the
locally supported standard FEM basis functions by basis functions formed by solutions of the $\mu$PDE with specific choices of parameters. By doing so, we typically obtain dense algebraic systems of considerably smaller size.

{The RBM is not a new approach and has been extensively studied in the literature for a wide range of applications, see, e.g.,~\cite{patera2006,quarteroni_rbm,hesthaven2015} and the references therein. 
  However, little attention has been paid to  applications in finance. Recent work on model reduction techniques in finance with a primary focus on POD methods
  include, e.g.,~\cite{sachs2008,sachs2014,sachs2014_1,pironneau12,red_bs15}. First results for RBM can be found in 
\cite{pironneau,pironneau11,pironneau12}. 
In ~\cite{MU14}, RBM was applied to more complex models, e.g., with parameter functions as an 
initial condition. While these references focus on the simplest case of European options,  American options, which are described by parabolic variational inequalites, are considered in~\cite{burkovska,burkovskaphd,Balajewicz2016}. 
In the case of parabolic variational inequalities, appropriately constructing the reduced basis spaces is much more challenging than for variational equalities.
To tackle this, POD-Angle-Greedy strategies~\cite{burkovska} or non-negative matrix factorization algorithms~\cite{Amsallem,Balajewicz2016} can be used. We also mention the relevant works on RBM for variational inequalities for the stationary case~\cite{HSW12,VeroyParisWorkshop} and instationary case in a space-time framework, e.g.,~\cite{UG13}. }

{Option pricing is associated with the calibration of non-observable parameters. This task can be formulated as a least-squares minimization problem, which finds a model parameter that minimizes the discrepancy
between the model and market option prices in the sense of least squares. These problems typically require require the numerical solution of the $\mu$PDE to be evaluated many times for different values of the parameter, which incurs to a high computational cost. Therefore, to accelerate the calibration routine while still providing accurate results, we aim to apply the RBM. The idea is to replace the  complex PDE model by a simpler surrogate model constructed by the RBM approaches, or alternatively by POD techniques, e.g.,~\cite{sachs2014,sachs2014_1}. Whereas calibration with the RBM in the linear case of European options has already been studied in the literature~\cite{pironneau,pironneau09}, the extension to American options, to the best of our knowledge, has not yet been addressed. We also mention the application of the RBM to PDE-constrained optimization problems~\cite{volkwein2010} using the POD method and recent work \cite{haasdonkopt} on the reduced basis method. }

{The calibration problem can be also
studied in the context of the theory of inverse problems. Most works are concerned with reconstructing the implied volatility surface for European options in the Black-Scholes model. This leads to an infinite-dimensional ill-posed inverse problem. We refer, e.g., to~\cite{bouchev}, where unique solvability and stability of the inverse problem are analyzed, and to~\cite{egger2005}, where an appropriate Tychonov regularization strategy is proposed to address the inherent ill-posedness of the problem.
Similarly, in~\cite{MR2137495,pir_calib}, an infinite-dimensional Tychonov regularized problem for American options is considered and the existence of solutions is derived together with their optimality conditions.  In the present case, due to the finite-dimensional
parameter space, we generally expect the corresponding inverse problem to
be well-posed, and so we do not use a regularization term.
%In the present finite-dimensional parameter domain and RB space, we do not encounter numerical stability problems and use no regularization term.
}

We note that for calibration one can also work directly in the stochastic
framework and compute the model prices, e.g., using a Monte-Carlo method by applying a backward regression scheme as in \cite{longstaff2001valuing} or by using different Monte-Carlo estimates as in e.g. \cite{broadie1997pricing,fu2001pricing,rogers2002monte}. Alternatively one can apply (binomial) tree methods as e.g. in \cite{crr_tree,rubinstein1994implied}. For European options, closed-form solutions can be used, or FFT techniques, see, e.g.,~\cite{mrazek2014optimization,Schoutens2004,MR2891932}. 
Fourier transform based pricing methods have been extended to price American
options, see for instance \cite{FangOosterlee2011} and \cite{Levendorskii2004}.
These methods are applicable when the Fourier transform of the modelling stochastic
process evaluated at fixed times is available in closed-form. In this article we
focus on a PDE-based approach that has a more general scope.

Whereas discretization techniques for PDEs and RDMs form a rather flexible and abstract framework, de-Americanization strategies are specially designed for the valuation of American options \cite{carr2010stock,deam2016}. These strategies first transform the price into a pseudo-European option price and secondly calibrate the European option by directly applying the computationally less expensive closed-form solution. In \cite{deam2016}, the de-Americanization techniques are studied numerically, which also reveals their limitations. 
{Both model reduction techniques, i.e., the RBM strategy and the de-Americanization strategy, are numerically studied and compared  for the Heston model~\cite{heston}. }

The rest of the  paper is structured as follows: The model problem and the calibration procedure is briefly introduced in \ref{ch6:sec:calib_heston}. In Section~\ref{ch6:sec:rb_calib}, we discuss a POD-Angle-Greedy RBM  for the variational inequality of an American option pricing problem based on the Heston model.  Section~\ref{ch6:sec:deam} is devoted to  the surrogate model obtained by the de-Americanization strategy (DAS). 
In Section~\ref{ch6:num_study} both RBM and DAS are investigated numerically.
A  comparative study of both techniques is presented in Section~\ref{ch6:sec:comparison}. Here, we use synthetic as well as  real market data sets to calibrate the unknown parameters.

\section{Model problem}
To calculate the price of an option (European or American) in a specific
model, suitable data must be provided. These data consist of the
current asset price $S_0>0$, the maturity time $T>0$, the strike price $K>0$, and the
set of input parameters, which are denoted by the vector $\parmu\in\p$, where $\p\subset\mathbb{R}^p$ is a parameter domain. 
Whereas
the first three components, $S_0$, $K$, $T$, are known and provided by the
market data, the input parameter vector $\parmu$ is not known a priori, except for the interest rate $r > 0$, and needs to be estimated from the
market. 
These types of problems are referred to as parameter identification  or calibration problems.
That is,
given a set of observations of market option prices, we are interested in
the parameter $\parmu$ that provides the best fit to the observed market data.

Generally, the market data is characterized by the actual  spot price $S_0$
 and the  market option prices
 $P^{\rm obs}_i=P^{\rm obs}(S_0,T_i,K_i)$ for different maturities
$T_i$ and for different strikes
$K_i$, $i=1,\dots,\Nopt$.
Mathematically, the calibration of option prices can be stated as a
least squares minimization problem: find $\parmu\in\popt\in\p$ that solves
\begin{equation}
 \min_{\parmu}J(\parmu),\quad\quad
J(\parmu):=\frac{1}{M}\sum_{i=1}^M|P^{\rm
obs}_i-P_i(\parmu)|^2,
\label{ch6:lsq}
\end{equation}
where 
$P_i(\parmu)=P(S_0,T_i,K_i;\parmu)$ are the model prices, which can for example be computed by solving the associated $\mu$PDE. More precisely, for European options, 
$P_i(\parmu)$ solves
\begin{subequations}
 \begin{align} 
\frac{\partial{P_i}(\parmu)}{\partial{\tau}}+\mathcal{L}(\parmu)P_i(\parmu)&=0,
&&\text{in}\quad [0,T)\times\mathbb{R}^n_+,\\
 P(T_i;\parmu)&=\mathcal{H}_i, &&\text{in}\quad\mathbb{R}^n_+,
\end{align}\label{ch6:calib_model_eu}\end{subequations}
where $\mathcal{H}_i$ is a pay-off functional and $\mathcal{H}_i=(K_i-S)_+$ for the put and 
$\mathcal{H}_i=(S-K_i)_+$ for the call options, $n=1,2$, and $\tau$ denotes the time to maturity $T:=\max(T_i)$, $i=1,\dots,M$. 
In the case of American put options, $P_i(\parmu)$ satisfies
\begin{subequations}
 \begin{align}
\frac{\partial{P_i}(\parmu)}{\partial{\tau}}+\mathcal{L}(\parmu)P_i(\parmu)&\leq
0,
&&\text{in}\quad[0,T)\times\mathbb{R}^n_+,\\
 P_i(\parmu)&\geq \mathcal{H}_i, &&\text{in}\quad[0,T)\times\mathbb{R}^n_+,\\
 \left(\frac{\partial
P_i(\parmu)}{\partial{\tau}}+\mathcal{L}(\parmu)P_i(\parmu)
\right)\left(P_i(\parmu)-\mathcal{H}_i\right)&=0,
&&\text{in}\quad[0,T)\times\mathbb{R}^n_+,\\
 P(T_i;\parmu)&=\mathcal{H}_i, &&\text{in}\quad \mathbb{R}^n_+,
\end{align}\label{ch6:calib_model_ao}\end{subequations}
with
$\mathcal{H}_i:=(K_i-S)_+$. Both problems are subject to
suitable boundary conditions. The operator $\mathcal{L}$ is a spatial (integro-) differential operator
and is defined by the model used to price the option, e.g., Black-Scholes~\cite{blackscholes}, CEV model, Heston~\cite{heston}.  
{Here, we restrict ourselves  to the Heston model,~\cite{MR2891932,mrazek2014optimization,Schoutens2004}, due to its ability to replicate the market behavior better than some other models. }

% H E S T O N
\section{The Heston model and its calibration}\label{ch6:sec:calib_heston}
We continue our discussion by introducing the Heston model~\cite{heston}, which is used to calibrate the option prices,
and and by describing the calibration procedure of this model.
\subsection{Model problem}
The Heston model is described by the following stock price~\eqref{ch2:mod_heston_s} and volatility~\eqref{ch2:mod_heston_v} 
dynamics,
\begin{subequations}
 \begin{align}
 dS&=\iota Sd\tau+\sigma SdW^1,\label{ch2:mod_heston_s}\\
 d\vol&=\kappa(\gamma-\vol)d\tau+\xi\sqrt{\vol}dW^2,\label{ch2:mod_heston_v}
\end{align}\label{ch2:mod_heston}\end{subequations}
The asset price $S:=\{S_\tau:\ \tau\geq 0\}$ exhibits geometric Brownian motion with Wiener process
$W^1$, drift $\iota$ 
and volatility $\sigma:=\sqrt{\vol}$. The stochastic instantaneous variance $\vol:=\{\vol_\tau:\ \tau\geq 0\}$
is driven by a mean-reverting square-root process (known as the Cox-Ingersoll-Ross (CIR) process) with long-run variance
$\gamma>0$, rate of mean reversion $\kappa>0$, and volatility of variance (also called the volatility of volatility) $\xi>0$.
The Wiener processes $W^1$ and $W^2$
are correlated by  $\rho\in[-1,1]$. 
Moreover, the so-called \textit{Feller condition} is  assumed, which states that the variance
process~\eqref{ch2:mod_heston_v} is strictly positive if the parameters satisfy
\begin{equation}
 2\kappa\gamma>\xi^2,\label{ch2:feller}
\end{equation}
see, e.g.,~\cite{janek2011fx}.
Unless otherwise stated, we focus on setting the parameters such that 
the Feller condition is always fulfilled. In terms of  It\^o calculus,  \eqref{ch2:mod_heston} can be reformulated as the partial differential equations
\eqref{ch6:calib_model_eu}, \eqref{ch6:calib_model_ao}, see, e.g.,~\cite{achdou}, where
$\parmu:=(\xi,\rho,\gamma,\kappa,r)$, and
the operator $\mathcal{L}(\parmu)$ is defined as follows
\begin{multline}
\mathcal{L}(\parmu)P(\parmu):=\frac{1}{2}\vol S^2\frac{\partial^2 P(\parmu)}{\partial S^2}+\xi
\vol\rho S\frac{\partial^2 P(\parmu)}{\partial\vol\partial
S}+\frac{1}{2}\xi^2\vol\frac{\partial^2 P(\parmu)}{\partial\vol^2}+rS\frac{\partial
P(\parmu)}{\partial S}\\+\kappa(\gamma-\vol)\frac{\partial
P(\parmu)}{\partial\vol}-rP(\parmu).
 \label{ch2:heston_operator}
\end{multline}
We note that the operator \eqref{ch2:heston_operator} is of a convection-diffusion reaction type with variable coefficients
that is degenerate for $S=0$, $\vol=0$.
To eliminate the degeneracy in the asset price variable, 
the standard practice is to 
perform a log-transformation of $S$. Introducing the new variable
$x:=\log(S/K)$, we denote by
$w(\tau,\vol,x,\parmu):=P(\tau,\vol,Ke^x;\parmu)$ the price in the log-transformed variable, i.e.,
$w(0,\vol,x;\parmu)=\chi(x):=(Ke^x-K)_+$ for a put and $\chi(x):=(K-Ke^x)_+$ for a call. Let  $L(\parmu)$ be the differential operator 
corresponding to  $\mathcal{L}(\parmu)$ in the log-transformed variable given by
\begin{subequations}\label{ch2:heston_operator_log_short2}
\begin{equation}
 L(\parmu)w(\parmu)=\nabla\cdot \Mat{A}(\parmu)\nabla w(\parmu)-\Mat{b}(\parmu)\cdot\nabla w(\parmu)-rw(\parmu),
 \label{ch2:heston_operator_log_short}
\end{equation}
with $\nabla :=\left(\frac{\partial}{\partial\vol},\frac{\partial}{\partial x}
\right)^T$, diffusion matrix $\Mat{A}(\parmu)$ and velocity vector $\Mat{b}(\parmu)$
\begin{eqnarray}\label{eq:heston_coef}
\Mat{A}(\parmu):=\frac{1}{2}\vol\begin{bmatrix}
\xi^2 &\rho\xi \\ \rho\xi &1
\end{bmatrix}, \ \ \ \ \ \ \ \ \ \ \ \ 
\Mat{b}(\parmu):=\begin{bmatrix}
-\kappa(\gamma-\vol)+\frac{1}{2}\xi^2 \\ -r+\frac{1}{2}\vol+\frac{1}{2}\xi\rho
\end{bmatrix}.
\end{eqnarray}\end{subequations}
The PDE is considered on a bounded
domain
$\Omega:=(\vol_{\min},\vol_{\max})\times(x_{\min},x_{\max})\subset\mathbb{R}^2$, 
$x_{\min}<0<x_{\max}$, $0<\vol_{\min}<\vol_{\max}$ with the Lipschitz continuous boundary $\partial\Omega=\Gamma_\DB\cup\Gamma_\NB$, where
$\Gamma_\DB$ corresponds to the non-trivial Dirichlet part of the boundary $\partial\Omega$, and $\Gamma_\NB$ stands for the Neumann part.
For European put options, we consider the boundary conditions as in~\cite{during},
\begin{subequations}
\begin{align}
\frac{\partial w}{\partial\vol}(t,\vol,x)=&0,&&\quad\text{on}\quad\Gamma_{\NB}:=\{(\vol,x)\in\Omega:\quad \vol\in\{\vol_{\min},\vol_{\max}\}\},\\
w(t,\vol,x)=&Ke^{-rt},&&\quad\text{on}\quad\Gamma_{\DB}^1:=\{(\vol,x)\in\Omega:\quad x=x_{\min}\},\\
w(t,\vol,x)=&0,&&\quad\text{on}\quad\Gamma_{\DB}^2:=\{(\vol,x)\in\Omega:\quad x=x_{\max}\}.
\end{align} \end{subequations}
For American put options, we specify the boundary conditions as proposed in~\cite{clarke,during},
\begin{subequations}
\begin{align}
w(t,\vol,x)=\chi(x),\quad\text{on}\quad\Gamma_{\DB}:=\Gamma_{\DB}^1\cup\Gamma_{\DB}^2, \quad
\frac{\partial w}{\partial\vol}(t,\vol,x)=0,\quad\text{on}\quad\Gamma_{\NB}.
\end{align} \end{subequations}
Next, we recast the problem in a variational form. We introduce the following functional spaces
\begin{equation}
 X=H^1(\Omega), \quad V:=\{v\in
X:\ v|_{\Gamma_\DB}=0\},\label{ch4:spaces}
\end{equation}
equipped with the norms $\|\cdot\|_X=\|\cdot\|_{H^1}$, $\|\cdot\|_V=|\cdot|_{H^1}$, which 
correspond to the $H^1(\Omega)$ norm and semi-norm respectively. Let $V^{\prime}$ be the dual space of $V$, and denote
by $\langle\cdot,\cdot\rangle_{\V\times V}$ the duality pairing of $V$ with $V^{\prime}$. We then define the bilinear 
form $a:V\times V\to\mathbb{R}$ corresponding to the Heston model, 
\begin{align}
  a(u,v;\parmu)&:=\int_{\Omega}\Mat{A}(\parmu)\nabla u\cdot\nabla
v+\int_{\Omega}\Mat{b}(\parmu)\cdot\nabla u v+\int_{\Omega}r u v.\label{ch4:bil_form_H}
\end{align}
Note that $\nu\geq\nu_{\min}>0$, $\rho\in(-1,1)$ and hence $\Mat{A}(\parmu)$
is positive definite on $\overline{\Omega}$. It follows from the admissible values of the parameters that for all $\parmu\in\p$ the 
bilinear form $a(\cdot,\cdot;\parmu)$ is continuous and satisfies a G\aa rding inequality on $V\times V$,~\cite{QV},
i.e., there exist constants
$0<\overline{\alpha}_a\leq\alpha_a(\parmu)$,
$0<\gamma_a(\parmu)\leq\overline{\gamma}_a<\infty$, $0\leq \lambda_a(\parmu)\leq\overline{\lambda}_a<\infty$, such that 
\begin{subequations}
\begin{align*}
 |a(u,v;\parmu)|&\leq \gamma_a(\parmu)\|u\|_V\|v\|_V \quad \forall u,v\in V,
&&&(\textit{continuity})\\
a(v,v;\parmu)&\geq\alpha_a(\parmu)\|v\|_V^2-\lambda_a(\parmu)\|v\|^2_\LL \quad \forall v\in V.
&&&(\textit{G\aa rding inequality})
\end{align*} \end{subequations}
%Moreover, under an appropriate conditions on the coefficients one can ensure the coercivity of $a(\cdot,\cdot;\parmu)$, see, e.g.,~\cite[Theorem 3.1]{winkler2001}.

For the semi-discretization in time, we 
use the $\theta$-scheme with $\theta=1/2$.
Let $t:=T-\tau$ and subdivide the time interval $[0,T]$ into $I$ subintervals of
equal length, $t^k:=k\Delta t$, $0<k\leq I$ with $\Delta t=T/I$. Defining
$w^k(\parmu):=w(t^k,\vol,x;\parmu)\in X$, we decompose it into  $w^k=u^k+u_L^k$, where
$u^k\in V$  and $u^k_L\in X$ is a continuous extension of the inhomogeneous Dirichlet data.
Furthermore, for all $\parmu\in\p$, $v\in V$, $\theta\in[0,1]$, $k\in\I$, where $\I:=\{0,\dots,I-1\}$, we define the linear functional
$f^{k+\theta}(\parmu)\in\V$ as follows
\begin{multline}
 f^{k+\theta}(v;\parmu):=
-\frac{1}{\Delta
t}\left(u^{k+1}_{L}(\parmu)-u^k_{L}(\parmu),v\right)_\LL\\
-a\left(\theta
u^{k+1}_{L}(\parmu)+(1-\theta)u_{L}^k(\parmu),v;\parmu\right).
\label{ch4:lin_form}
\end{multline}
Using a test function $v\in V$ in \eqref{ch6:calib_model_eu} yields the following variational formulation
\begin{align} 
  \frac{1}{\Delta t}\left(u^{k+1}-u^k,v\right)_\LL+a(u^{k+\theta},v;\parmu)
  =f^{k+\theta}(v;\parmu), \quad v\in V,
  \label{var_form_main}
 \end{align}
where $u^{k+\theta}:=\theta u^{k+1}+(1-\theta)u^k$. 

In terms of a  Lagrange multiplier,  the constrained problem~\eqref{ch6:calib_model_ao} can be weakly written in saddle point form, see, e.g.,~\cite{kikuchi}.
Define $W=V^{\prime}$ and $M\subset W$ as a dual cone by
\begin{equation}
 M:=\{\eta\in W:\quad b(\eta,v)\geq 0, \ v\in V, \ v\geq 0\},
\end{equation}
where  $b(\eta,v):=\langle\eta,v\rangle_{\V\times V}$, for all $\eta\in W$, $v\in V$.  
For $\parmu\in\p$ and $k\in\I$, we introduce the functional
$g^{k+1}\in W^{\prime}$, 
\begin{equation*}
 \g^{k+1}(\eta;\parmu):=b(\eta,\pf)-b(\eta,u_L^{k+1}
(\parmu)).
\end{equation*}
Then we arrive at the following variational saddle point formulation: For $\parmu\in\p$, $k\in\I$, $\theta\in[0,1]$, 
find $(u^{k+1}(\parmu),\lambda^{k+1}(\parmu))\in V\times M$, satisfying for all $\eta\in M, v\in V$,
\begin{align}
  \frac{1}{\Delta t}\left(u^{k+1}-u^k,v\right)_\LL+a(u^{k+\theta},v;\parmu)
  - b(\lambda^{k+1},v)  =  f^{k+\theta}(v;\parmu),\\
  b(\eta-\lambda^{k+1},u^{k+1})  \geq
  g^{k+1}(\eta-\lambda^{k+1};\parmu). 
  \label{var_form_vi}
\end{align}

Note that for all $\parmu\in\p$, $k\in\I$, $f^{k+\theta}\in\V$ and the bilinear form $a(\cdot,\cdot;\parmu)$ is continuous and 
satisfies the G\aa rding inequality. Thus, for  a small enough time step $\Delta t<(1/\theta\lambda_a(\parmu))$,
by a generalized Lax-Milgram argument, the problem~\eqref{var_form_main}
admits a unique solution $u(\parmu)\in V$,~\cite{achdou,QV}. 
Moreover, the bilinear from $b(\cdot,\cdot)$ is inf-sup stable
on $W\times V$ and  unique solvability of ~\eqref{var_form_vi} is given, see, e.g.,~\cite[Theorem 2.1]{brezzi_raviart}.

\subsection{Calibration}\label{ch6:sec:model} 
Given $S_0$ and the observations $P_i^{\rm obs}$ at $(T_i,K_i)$,
$i=1,\dots,M$, we need to compute the corresponding prices in the Heston
model $P_i(\vol_0,S_0,T_i,K_i)$, where $\vol_0\in\mathbb{R}_+$ is the initial
volatility. However, the  value of $\vol_0$ is not observable and needs to be determined
together with the parameters $\xi,\rho,\gamma,\kappa$. We collect all parameters to be
identified into the single
vector
$\muopt=(\muopti_1,\dots,\muopti_5)=(\xi,\rho,\gamma,\kappa,\vol_0)\in\popt$,
where
\begin{equation}
 \popt:=\left\{\muopt\in\mathbb{R}^5:
\muopti_{\min,i}\leq\muopti_i\leq\muopti_{\max,i},\quad i=1,\dots,5\right\}.
\end{equation}
We replace the PDE constraints with the corresponding weak
discrete problem in a log-transformed variable. The polygonal domain $\Omega\subset\mathbb{R}^2$ is resolved by 
a triangulation $\mathcal{T}_\N$ of $\Omega$, consisting of
$J$ simplices $T_\N^j$, $1\leq j\leq J$, such
that $\overline{\Omega}=\cup_{T_\N\in\mathcal{T}_\N}\overline{T}_\N$. We use standard
conforming nodal first-order finite element approximation spaces $X_\N\subset X$, $V_\N\subset V$, where
$X_\N:=\{v \in X:\quad v_{|T_\N^j}\in\mathbb{P}^1(T_\N^j),
1\leq j\leq J \}$, and $\mathbb{P}^1$ is a space of linear polynomials with degree at most one, and $V_\N=X_\N\cap V$, $\dim(X_\N)=\N_X$,
$\dim(V_\N)=\N_V$. 
To discretize the dual space
$W$, we use discontinuous piecewise linear biorthogonal basis
functions defined on the
same mesh as the basis functions of $V_\N$, \cite{Woh00a}.
That
is, $W_\N:=\spann\{\chi_q,\ q=1,\dots, \N\}$, $\dim(W_\N)=\N_W = \N_V$, where $\dual_q$ satisfy
a local biorthogonality relation:
\begin{equation*}
 \int_{T_\N^j}\chi_q\phi_p=\delta_{pq}\int_{T_\N^j}\phi_p\geq 0, \quad\phi_p\in
V_\N, \quad p,q=1,\dots,\N_V.
\end{equation*}
The discrete cone $M_\N\subset W_\N$ is given as
\begin{equation}
 M_{\N}=\spann_+\{\chi_q\}_{q=1}^{\N_W}:=\Big\{\eta\in W_\N:\quad
\eta=\sum_{q=1}^{\N_W} \alpha_q\chi_q ,
\quad \alpha_q\geq 0\Big\}
\end{equation}
and $(V_\N, W_\N)$ form a uniformly inf-sup stable pairing.

For a given $\parmu\in \p$, $k\in\I$, we approximate
$w^k(\parmu)$ by $w_\N^k(\parmu)\in X_\N$, $w_\N^k(\parmu):=u_\N^k(\parmu)+u_{L\N}^k(\parmu)$ with
$u_\N^k(\parmu)\in V_\N$ and $u_{L\N}^k(\parmu)\in X_\N$
and  $\lambda^{k+1}(\parmu) $ by $\lambda_\N^{k+1}(\parmu)\in M_\N$, $k\in\I$.
In the reduced basis literature it is common to call $u_\N^k(\parmu)$ and $\lambda_\N^k(\parmu)$ the  detailed solutions.

Then, for the European put options $u_\N^{k+1}(\parmu)\in V_\N$, $\parmu\in\p$, $k\in\I$,
solves the following discrete problem 
 \begin{align*} 
  \mathbb{E}_\N^{\rm Eu}(\parmu) = \begin{cases}
  \frac{1}{\Delta t}\left(u_\N^{k+1}-u_\N^k,v\right)_\LL+a(u_\N^{k+\theta},v;\parmu)
  =f^{k+\theta}(v;\parmu), \\ 
  (u_\N^0 - u^0, v)_V = 0,\qquad v\in V_\N.
  \end{cases}
 \end{align*}
The solution pair $(u_\N^{k+1}(\parmu),\lambda_\N^{k+1}(\parmu))\in
V_\N\times M_\N$, $\parmu\in\p$, $k\in\I$, for the American option problem in turn satisfies the following detailed saddle point problem
\begin{align*}
  \mathbb{E}_\N^{\rm Am}(\parmu) =
  \begin{cases}
  \frac{1}{\Delta t}\left(u_\N^{k+1}-u^k_\N,v\right)_\LL+a(u_\N^{k+\theta},v;\parmu)
  - b(\lambda_\N^{k+1},v)  =  f^{k+\theta}(v;\parmu),\\
  b(\eta-\lambda_\N^{k+1},u_\N^{k+1})  \geq
  g^{k+1}(\eta-\lambda_\N^{k+1};\parmu),\qquad \eta\in M_\N, v\in V_\N,\\
  (u_\N^0 - u^0, v)_V = 0.
  \end{cases} 
\end{align*}
Both problems have a unique discrete solution,
and for $i=1,\dots,M$  we denote the approximate model price as
$$P_i^{\N,\rm
s}(\muopt):=w_\N^{k_i}(\log(S_0/K_i),\vol_0;\parmu),$$ 
which is obtained from the solutions of $\mathbb{E}^{s}_\N$, $s \in \{\rm Eu, \rm
Am \}$. By abuse of notation, we often omit the index $s$ in the
notation of $P_i^{\N,s}$, if it is clear from context.  
Then the minimization
problem~\eqref{ch6:lsq} can be written in the following form 
\begin{equation}
\begin{aligned}
 &\min_{\muopt\in\popt}J_\N(\muopt):=\frac{1}{M}\sum_{i=1}^M|P^{\rm
obs}_i-P^\N_i(\muopt)|^2.\\
% &\text{subject to}\quad \mathbb{E}_\N(\parmu),
\end{aligned}\label{ch6:minprob}
\end{equation}

\section{Reduced basis method (RBM)}~\label{ch6:sec:rb_calib}
The high-fidelity discrete problem
$\mathbb{E}_\N$, in general, is computationally expensive for large $\N_V$, and significantly slows down the calibration procedure. To reduce the complexity, we apply the
reduced basis method and
substitute the detailed model $\mathbb{E}_\N$ with the reduced-order surrogate
model $\mathbb{E}_N$. 
Using a suitable
basis generation procedure, as discussed below, we construct the 
low-dimensional reduced
spaces $V_N\subset V_\N$ for European options and $V_N\subset V_\N$,
$W_N\subset W_\N$, $M_N\subset
M_\N$ for American options with $\dim(V_N)\ll\dim(V_\N)$,
$\dim(W_N)\ll\dim(W_\N)$. For a given $\parmu\in\p$, we approximate $u_\N^{k+1}\in V_\N$ by $u_N^{k+1}(\parmu)\in V_N$, 
$\lambda_\N^{k+1}(\parmu)\in M_\N$ by $\lambda_N^{k+1}(\parmu)\in M_N$, $k\in\I$.
The
reduced surrogate models for pricing European and American options are defined as
follows
\begin{align*} 
  \mathbb{E}_N^{\rm Eu}(\parmu)& = \begin{cases} \frac{1}{\Delta
t}\left(u_N^{k+1}-u_N^k,v\right)_\LL+a(u_N^{k+\theta},v;\parmu)=f^{
k+\theta}(v;\parmu), \\ 
(u_N^0-u_\N^0,v)_V=0,\qquad v\in V_N. \end{cases}\\
  \mathbb{E}_N^{\rm Am}(\parmu)&=
  \begin{cases}
\frac{1}{\Delta
t}\left(u_N^{k+1}-u^k_N,v\right)_\LL+a(u_N^{k+\theta},v;\parmu)
-b(\lambda_N^{k+1},v)  =  f^{k+\theta}(v;\parmu),\\
b(\eta-\lambda_N^{k+1},u_N^{k+1})  \geq
g^{k+1}(\eta-\lambda_N^{k+1};\parmu),\qquad \eta\in M_N,\ v\in V_N,\\
(u_N^0-u_\N^0,v)_V=0.\end{cases} 
\end{align*}
Additionally, we require that the reduced spaces $V_N$, $W_N$ are constructed such that the bilinear form $b(\cdot,\cdot)$ is uniformly inf-sup stable on $W_N\times V_N$ with respect to $N$. Thus the well-posedness of the reduced problems $\E_N^{\rm Eu}(\parmu)$ and $\E_N^{\rm Am}(\parmu)$ is given, see also~\cite{burkovska,HSW12}. 

We denote the reduced basis approximation of the model price by
$$P_i^{N,s}(\muopt):=w_N^{k_i}(\log(S_0/K_i),\vol_0;\parmu),$$ 
where 
$w_N^{k_i}(\parmu):=u_N^{k_i}(\parmu)+u_{L\N}^{k_i}(\parmu)$ and $u_N^{k_i}(\parmu)$ is a solution of
$\mathbb{E}_\N:=\mathbb{E}^{s}_\N$, $s=\{\rm Eu, Am\}$ for $i=1,\dots,M$. Then
we approximate $J_\N(\muopt)\approx J_N(\muopt)$, and obtain the following
minimization problem for the reduced model
\begin{equation}
\begin{aligned}
 &\min_{\muopt\in\popt}J_N(\muopt):=\frac{1}{M}\sum_{i=1}^M|P^{\rm
obs}_i-P^N_i(\muopt)|^2.\\
% &\text{subject to}\quad \mathbb{E}_N(\parmu),
\end{aligned}\label{ch6:minprob_rb}
\end{equation}

\begin{remark}
Note that the interest rate $r$ is not determined by a calibration procedure and is
fixed beforehand. In our case, the market data will be a single stock in the U.S., and for an approximation of the risk-free rate we use the rates of the U.S. Department of the Treasury. Hence, to construct the reduced basis
spaces, we only need to consider the variation of four parameters
$\parmu=(\xi,\rho,\gamma,\kappa)\in\p\subset\mathbb{R}^4$. However, this choice
is restrictive and for new market data we would need to construct a new
reduced basis set. Therefore, to conserve generality in our approach, we
consider the variation of all parameters, $\parmu=(\xi,\rho,\gamma,\kappa,r)$, and
consequently the constructed reduced basis will be entirely 
market-independent.
\end{remark}

Numerous approaches exist for the construction of the reduced basis approximation spaces. Their common goal is to exploit the parameter dependence of the problem and to incorporate this information into the construction of the reduced bases. Typically, this is done by applying an iterative greedy strategy to a set of snapshots, i.e., solutions computed for different parameter values. For linear parabolic problems, a popular choice is a combination of the greedy strategy for parameter selection and a proper orthogonal decomposition (POD) in time resulting in a so-called POD-Greedy algorithm~\cite{Ha13,haasdonk08}.
For parabolic variational inequalities, the construction is more challenging due to the requirement of uniform inf-sup stability. For stationary variational inequalities, a greedy sampling is commonly used,~\cite{HSW12,VeroyParisWorkshop}, while for time-dependent problems a POD-Angle-Greedy~\cite{burkovska} and a POD-NNMF\footnote{NNMF refers to a non-negative matrix factorization procedure,~\cite{Lee_Seung}.}~\cite{Amsallem} have been considered in the literature. In the present work, we follow the idea of the POD-Greedy  algorithm for European options and POD-Angle-Greedy algorithm for American options. 

Since the European option problem can be considered as a particular case of the American option problem formulation,  we focus on the description of the basis construction for the latter problem and comment only on the differences. 

Consider a finite subset $\p_N:=\{\parmu_1,\dots,\parmu_N\}$, $\parmu_i\neq\parmu_j$, $\forall i\neq j$, $N\in\mathbb{N}$ and define the reduced spaces $V_N:=\spann\{\Psi_N\}$ and $W_N:=\spann\{\Xi_N\}$, where the primal $\Psi_N:=\{\psi_1,\dots,\psi_{N_V}\}\subset V_\N$, and 
the dual $\Xi_N:=\{\xi_1,\dots,\xi_{N_W}\}\subset W_\N$ reduced bases are constructed from the large set of snapshot solutions $u_\N^k(\parmu_i)$ and $\lambda_\N^{k+1}(\parmu_i)$, $k\in\I$, $i=1,\dots,N$. The reduced cone is defined as $M_N:=\spann_+\{\xi_j\}_{j=1}^{N_W}:=\left\{\sum_{j=1}^{N_W}\alpha_j\xi_j, \ \alpha_j\geq 0\right\}$. By construction $\xi_j\in M_\N$ and thus $M_N\subset M_\N$. 

The approach we follow to construct $\Psi_N$ and $\Xi_N$ is presented in Algorithm~\ref{ch5:alg:angle-greedy}. We investigate the parameter domain $\p$, that is replaced by a finite    set $\p_{\train}\subset\p$, by a greedy search (Step~\ref{ch5:alg:angle-greedy5}--Step~\ref{ch5:alg:angle-greedy13}). In the greedy loop, we identify a ``worst'' parameter $\parmu_N$, i.e.,  a parameter which leads to the worst RB approximation, and add it to the training set. This selection requires an error measure $E_N(\parmu_N)$, which can be chosen, e.g., as the true error between the detailed solution and the reduced basis approximation or an a posteriori error bound~\cite{burkovska,grepl2005reduced,haasdonk08}. 
The availability and efficient computation of the latter choice makes it more attractive from a computational point of view. 
Next, for the selected parameter, we compute primal and dual bases and repeat this process $N_{\max}$ times or until the desired tolerance $\varepsilon_{\tol}$ of the stopping criterion is reached. 
For the primal reduced space construction, we apply the standard POD-Greedy procedure (Step~\ref{ch5:alg:angle-greedy11}), where the difference between the worst resolved trajectory $u^k_\N(\parmu_N)$, $k\in\I$, and its orthogonal projection onto the current RB space $\Pi_{V_{N-1}}u_\N^k(\parmu_N)$ is compressed to the first dominant POD mode:
$$ POD_1\left(\{v^k\}_{k \in\I_0}\right):=\argmin_{\|z\|_V=1}\sum_{k\in \I_0}\|v^k-(v^k,z)_Vz\|_V^2.$$
 \begin{algorithm}[H]
\caption{POD-Angle-Greedy Algorithm}
\label{ch5:alg:angle-greedy}
\begin{algorithmic}[1]
\Require Maximum number of iterations $N_{\max}>0$, training sample set
$\p_{\train}\subset\p$, target tolerance $\varepsilon_{\tol}$
\Ensure RB bases $\Psi_N$, $\Xi_N$ and RB spaces $V_N$, $W_N$
\State arbitrarily choose $\parmu_0\in\p_{\train}$ and
$k^{\prime}\in\I_0:=\I\cup\{I\}$\label{ch5:alg:angle-greedy0} 
\State compute $\{u_\N^{k}(\parmu_0)\}_{k\in\I_0}$,
$\{\lambda_\N^{k+1}(\parmu)\}_{k\in\I}$
\State set
$\xi_0={\lambda_\N^{k^{\prime}}(\parmu_0)}/{\|\lambda_\N^{k^{
\prime } } (\parmu_0)\|_W}$, $\Xi_0=\{\xi_0\}$,
$W_0=\spann\{\Xi_0\}$
\State set
 $\Psi_0={\rm
orthonormalize}\left\{u_\N^{k^{\prime}}(\parmu_0),\B\xi_0\right\}$,
$V_0=\spann\{\Psi_0\}$\label{ch5:alg:angle-greedy4} 
\For{$N=1,\dots,N_{\max}$}\label{ch5:alg:angle-greedy5} 
    \State $
\parmu_{N}=\arg\max_{\parmu\in\p_{\train}}{\E_{N-1}
(\parmu) } $\label{ch5:alg:angle-greedy6}
 \If{$\varepsilon^{\train}_{N}<\varepsilon_{\tol}$}
 \label{ch5:alg:angle-greedyep}
	\Return
    \EndIf
\State $k_N=\arg\max_{k\in\I}\left( \measuredangle \left(
\lambda_\N^{k+1}(\parmu_N) ,
W_{N-1}\right)\right)$\label{ch5:alg:angle-greedy9}
\State
$\xi_N={\lambda_\N^{k_N}(\parmu_N)}/{\|\lambda_\N^{k_{N} }
(\parmu_N)\|_W}$, $\Xi_N=\Xi_{N-1}\cup\{\xi_N\}$,
$W_N=\spann\{\Xi_{N}\}$\label{ch5:alg:angle-greedy10}
    \State $\psi_{N}=POD_1\left(\left\{u_\N^k(\parmu_{N})-\Pi_{
V_{N-1}}\left(u_\N^k(\parmu_{N})\right)\right\}_{k\in\I_0}\right)$
\label{ch5:alg:angle-greedy11}
\State
$\Psi_{N}={\rm
orthonormalize}\left\{\Psi_{N-1}\cup\{\psi_N,\B\xi_N\}\right\}$,
$V_{N}=\spann\{\Psi_{N}\}$\label{ch5:alg:angle-greedy12}
\EndFor \label{ch5:alg:angle-greedy13} 
\end{algorithmic}
\end{algorithm}

To construct the dual RB space, the vectors that maximize the volume of the resulting cone, i.e., vectors showing the largest deviation from the current RB space (Step~\ref{ch5:alg:angle-greedy9}, are selected. We denote $\measuredangle(\eta,Y):=\arccos\left(\|\Pi_Y\eta\|_W/\|\eta\|_W\right)$ the angle between the vector $\eta\in W$ and the linear space $Y\subset W$, where $\Pi_Y\eta$ is an orthogonal projection of $\eta$ onto $Y$. 

Additionally, to form a uniformly stable pair of the reduced spaces $V_N$, $W_N$,  we enrich the primal space $V_N$ by the ``supremizing'' $T\xi_N$,~\cite{rovas,rozza_2007}, where $T:W_\N\to V_\N$ is a ``supremizing'' operator, defined as a solution of $(T\xi_N,v)_V:=b(\xi_N,v)$, for all $v\in V_\N$. It is easy to see that this inclusion ensures the inf-sup stability of $b(\cdot,\cdot)$ on $W_N\times V_N$, see, e.g.,~\cite{HSW12,rozza_2007,HSW12}. 

For the case of European options no dual space is required and thus the steps  in Algorithm~\ref{ch5:alg:angle-greedy} involving the Lagrange multiplier space are omitted, resulting in a standard POD-Greedy algorithm.

A computational speed-up of this method is achieved by a so-called offline/online procedure, which relies on the assumption that the problem has an affine dependence on its parameters. That is, for every $\parmu\in\p$ the bilinear and linear forms are separable, i.e., there exist $\Theta_q^a:\p\to\mathbb{R}$, $q=1,\dots,Q_a$, such that $a(\cdot,\cdot;\parmu):=\sum_{q=1}^{Q_a}\Theta_q^a(\parmu)a_q(\cdot,\cdot)$, where $a_q:V\times V\to\mathbb{R}$ are parameter-independent. The same arguments apply to $f^{k+\theta}(\cdot;\parmu)$, $g^k(\cdot;\parmu)$ and $u_\N^0(\parmu)$. Then an offline routine requires the evaluation of all parameter-independent quantities. Usually this procedure is computationally cost-intense, however it is performed only once. In turn, the online procedure is computationally cheap and involves assembling the parameter-dependent components and solving the reduced system. This stage is executed multiple times for each new parameter value $\parmu\in\p$; see, e.g.,~\cite{hesthaven2015,quarteroni_rbm} for more details.

\section{De-Americanization strategy (DAS)}\label{ch6:sec:deam}
 In contrast to the reduced basis framework, which is based on the underlying pricing PDE,  the DAS transforms the observed American market prices into pseudo-European prices prior to the calibration.
That is, given an input data of American put options, we consider
the minimization problem~\eqref{ch6:minprob} with $J_\N(\muopt)\approx
\deam{J}_\N(\muopt)$,
\begin{equation}
\begin{aligned}
 &\min_{\muopt\in\popt}\deam{J}_\N(\muopt):=\frac{1}{M}\sum_{i=1}^M|\deam{
P}^ { \rm
obs}_i-P^\N_i(\muopt)|^2,\\
 %&\text{subject to}\quad \mathbb{E}^{\rm Eu}_\N(\parmu),
\end{aligned}\label{ch6:minprob_deam}
\end{equation}
where the prices $\deam{P}^{\rm obs}_j$ are the pseudo-European put
prices. These are obtained by perturbing the American put prices
$P_j^{\rm obs}$, i.e., $\deam{P}^{\rm obs}_j:=\tdeam(P_j^{\rm obs})$, where
$\tdeam:\mathbb{R}^M\to\mathbb{R}^M$, and the
corresponding model prices are the European put prices,
$P_j^\N(\muopt)=P_j^{\N,\rm
Eu}(\muopt)$. 
As in \cite{deam2016}, we use the binomial tree method, see \cite{crr_tree}, to transform American option prices into pseudo-European option prices. In a nutshell, once the observed market prices $P_j^{\rm obs}, j=1,\ldots,M$, have been collected, for each single stock option an individual binomial tree is calibrated to match this option price $P_j^{\rm obs}$. The resulting tree is used to price the associated  so-called pseudo-European option with the same strike and maturity. A detailed description of the de-Americanization method is given in \cite[Algorithm 1]{deam2016}.

The advantage of this method is that the complexity of the non-linear model for pricing American options can be reduced to that of the linear model for pricing European options, allowing closed-form solutions or Fourier techniques to be exploited.

In the following, we denote the prices obtained by the semi-closed-form solutions by $P_i^{\rm CF}$, see~\cite{heston}.
This leads to the following minimization problem
\begin{equation}
\begin{aligned}
 &\min_{\muopt\in\popt}\deam{J}_{\rm CF}(\muopt):=\frac{1}{M}\sum_{i=1}^M|\deam{
P}^ { \rm
obs}_i-P^{\rm CF}_i(\muopt)|^2.
 %&\text{subject to}\quad \mathbb{E}^{\rm Eu}_\N(\parmu),
\end{aligned}\label{ch6:minprob_deam_cf}
\end{equation}
We note that from the computational point of view, the DAS is
very attractive, particularly in combination with the closed-form solutions.
However, for each set of observations, an
additional pre-processing time to transform the American data into European
one is required, and, as we will see later, this step can dominate  the
computational
cost of the entire calibration routine. 

One could also
consider a combination of the RBM with the de-Americani\-zation strategy, i.e.,
applying the RBM to approximate $E_\N^{\rm Eu}$ by $E_N^{\rm Eu}$. The
corresponding minimization problem can be stated as follows
\begin{equation}
\begin{aligned}
 &\min_{\muopt\in\popt}\deam{J}_N(\muopt):=\frac{1}{M}\sum_{i=1}^M|\deam{
P}^ { \rm
obs}_i-P^N_i(\muopt)|^2,\\
 %&\text{subject to}\quad \mathbb{E}^{\rm Eu}_\N(\parmu),
\end{aligned}\label{ch6:minprob_deam_rb}
\end{equation}
with $P_i^N(\muopt)=P_i^{N,\rm Eu}(\muopt)$. {We note that due to the presence of the box constraints this finite-dimensional minimization problem admits a solution, see, e.g.,~\cite{vexler}.} %\todo{is this really correct ??}

\section{Numerical results}\label{ch6:num_study}

The problems under consideration belong to the class of
finite-dimensional
optimization
problems with box constraints and thus suitable numerical algorithms have to be applied.
For the special case of the calibration with European options, the most popular optimization algorithms are the gradient-based optimization methods; see,
e.g.,~\cite{achdou,sachs2014} and the references therein. 
By contrast, for American options, the situation is more involved  due to non-differentiability, see, e.g.,
\cite{ito2000,hintermuller,schiela}.
Here, we use the MATLAB Optimization
Toolbox and apply the built-in ``black-box''
optimization solver {\it lsqnonlin} or {\it fmincon}, in which the gradients
are approximated by finite differences. 

For the numerical experiments we set
$T=2$, $I=250$, $\Delta t=T/I=0.008$, $\theta=1/2$. The computational domain
$\Omega=(\vol_{\min},\vol_{\max})\times(x_{\min},x_{\max})=(10^{-5},
3)\times(-5,5)$ is resolved by a triangulation with $\N_X=4753$ nodes. For
$\parmu:=(\xi,\rho,\gamma,\kappa,r)\in\p\subset\mathbb{R}^5$ and
$\muopt:=(\xi,\rho,\gamma,\kappa,\vol_0)\in\p^{\rm opt}\subset\mathbb{R}^5$, we
define 
\begin{align}
\p&\equiv[0.1,0.9]\times[-0.95,0.95]\times[0.01,0.5]\times[0.1,
5]\times[0.0001,0.8],\\
\p^{\rm opt}&\equiv[0.1,0.9]\times[-0.95,0.3]\times[0.01,0.5]\times[0.1,
5]\times[10^{-5},1].
\end{align}
{Unless otherwise stated, the calibration routine is performed with with ${\it lsqnonlin}$, which uses
a Trust-Region-Reflective algorithm, and the stopping criterion is set as
$J(\muopt)-J(\muopt^{\star})\leq10^{-12}$, $\|\muopt-\muopt^{\star}\|_2\leq10^{-5}$, where $\muopt^{\star}$ is a locally optimal solution.}
%\todo{where to you have $\muopt^{\star} $ from in pratice you do not have it otherwise there would be no need for doing numerics}

\subsection{Calibration based on RBM}\label{ch6:numerics:rbm}
We consider a training set
$\p_{\train}$ composed
of uniformly distributed points in $\p$ with
$|\p_{\train}|=1024$ for the European put and $|\p_{\train}|=3125$ for the
American put options. The bases are generated by
the POD-Greedy and
POD-Angle-Greedy algorithms for the European and the American options, respectively.
The reduced systems have dimension $N_{\max}=100$ for European put and $N_{\max}=125$ for American
put options.
Firstly, we consider the quality of the calibration in terms  of the RBM applied to a synthetic data set with $r=5\%$:

\begin{equation}
\begin{alignedat}{2}
S_0&=1,\\
T_1&=\frac{1}{6},\qquad &&K_1=\{0.95, 0.975, 1, 1.025, 1.05\},\\
T_2&=\frac{1}{2},\qquad &&K_2=K_1\cup\{0.9,0.925,1.075,1.1\},\\
T_3&=\frac{3}{4},\qquad &&K_3=K_2\cup\{0.85,0.875,1.125,1.15\},\\
T_4&=1,\qquad &&K_4=K_3\cup\{0.8,0.825,1.175,1.2\},\\
T_5&=2,\qquad &&K_5=K_4\cup \{0.75, 0.775,1.225,1.25\}.
\end{alignedat}
\label{eq:syndatasetting}
\end{equation}
For each pair $(T_i,K_i)_{i=1,\dots,5}$, we
generate two artificial sets of observations $P^{\rm obs}$ consisting of 65 European and American put options at 
$\muopt=(0.7,-0.8,0.3,1.4,0.3)$.
That is, we solve the detailed problems~\eqref{ch6:minprob} with $\mathbb{E}_\N^{\rm Am}$ and $\mathbb{E}_\N^{\rm Eu}$
for the parameter $\parmu=(0.7,-0.8,0.3,1.4,0.05)$ and $K=1$ and interpolate the corresponding solution
$K_iu_\N^{k_i}(\vol,x;\parmu)$ at  $\vol^{\star}=\vol_0$, $x_i^{\star}=\log(S_0/K_i)$. 

We perform the optimization routine with the reduced surrogate
model~\eqref{ch6:minprob_rb} and the high-fidelity detailed problem~\eqref{ch6:minprob}. In both cases, we use the same initial guess  $\muopt_{\rm in}$ for the optimization algorithm. The results of the calibration
for two different data sets of American and European options are presented in
Table~\ref{calib_table_smt}. We observe that, using the detailed models
$\mathbb{E}_\N^{s}$, $\{s=\rm Am, Eu\}$, we can basically recover the exact
parameters $\muopt_{\rm ex}$. The reduced
surrogate models provide  still accurate enough results but are computationally much less expensive.

\begin{table}[h]
   \centering
\begin{tabular}{ccccccccc}                         
\toprule                                                
Method&$\mathbb{E}(\parmu)$ & $\muopt$& $\xi$ &$\rho$ &$\gamma$ &$\kappa$& $\vol_0$& $\|\muopt_{\rm ex}-\muopt^{\star}\|_2$ \\ 
\midrule
 &&$\muopt_{\rm ex}$ &0.700&   -0.800&    0.300&    1.400&    0.300& \\
 &&$\muopt_{\rm in}$ &0.601 &-0.682& 0.487 &2.020 &0.496& \\        
 \midrule
  $J_\N(\muopt)$&$\mathbb{E}_\N^{\rm Am}$&$\muopt^{\star}$ &0.700&   -0.800&    0.300&    1.399&    0.300&2.14e-5 \\
 $J_N(\muopt)$&$\mathbb{E}_N^{\rm Am}$&$\muopt^{\star}$ &0.694 &-0.831& 0.298 &1.447 &0.303&5.62e-2 \\
 $J_\N(\muopt)$&$\mathbb{E}_\N^{\rm Eu}$&$\muopt^{\star}$ &0.700&   -0.800&    0.300&    1.399&    0.300& 2.05e-5\\
 $J_N(\muopt)$&$\mathbb{E}_N^{\rm Eu}$&$\muopt^{\star}$ &0.616&   -0.886&    0.293&    1.306&    0.300&1.52e-1 \\
\bottomrule
\end{tabular} 
\caption{Calibration of the synthetic data set of American and European put options in the Heston model using detailed and reduced problems.}
\label{calib_table_smt}
\end{table}

\begin{figure}[ht]
\centering
 \includegraphics[width=.45\linewidth]
 {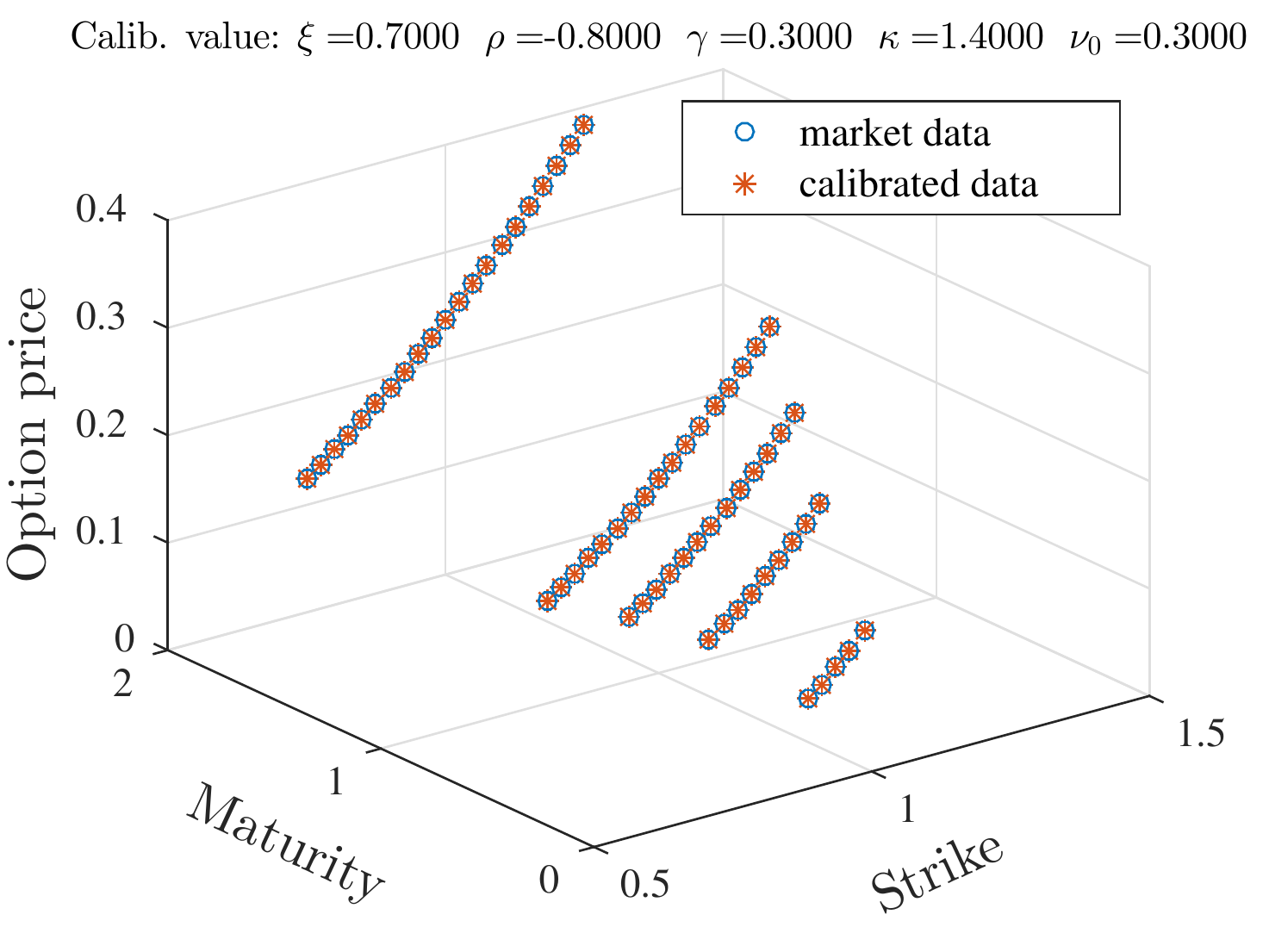}
 \includegraphics[width=.45\linewidth]
                 {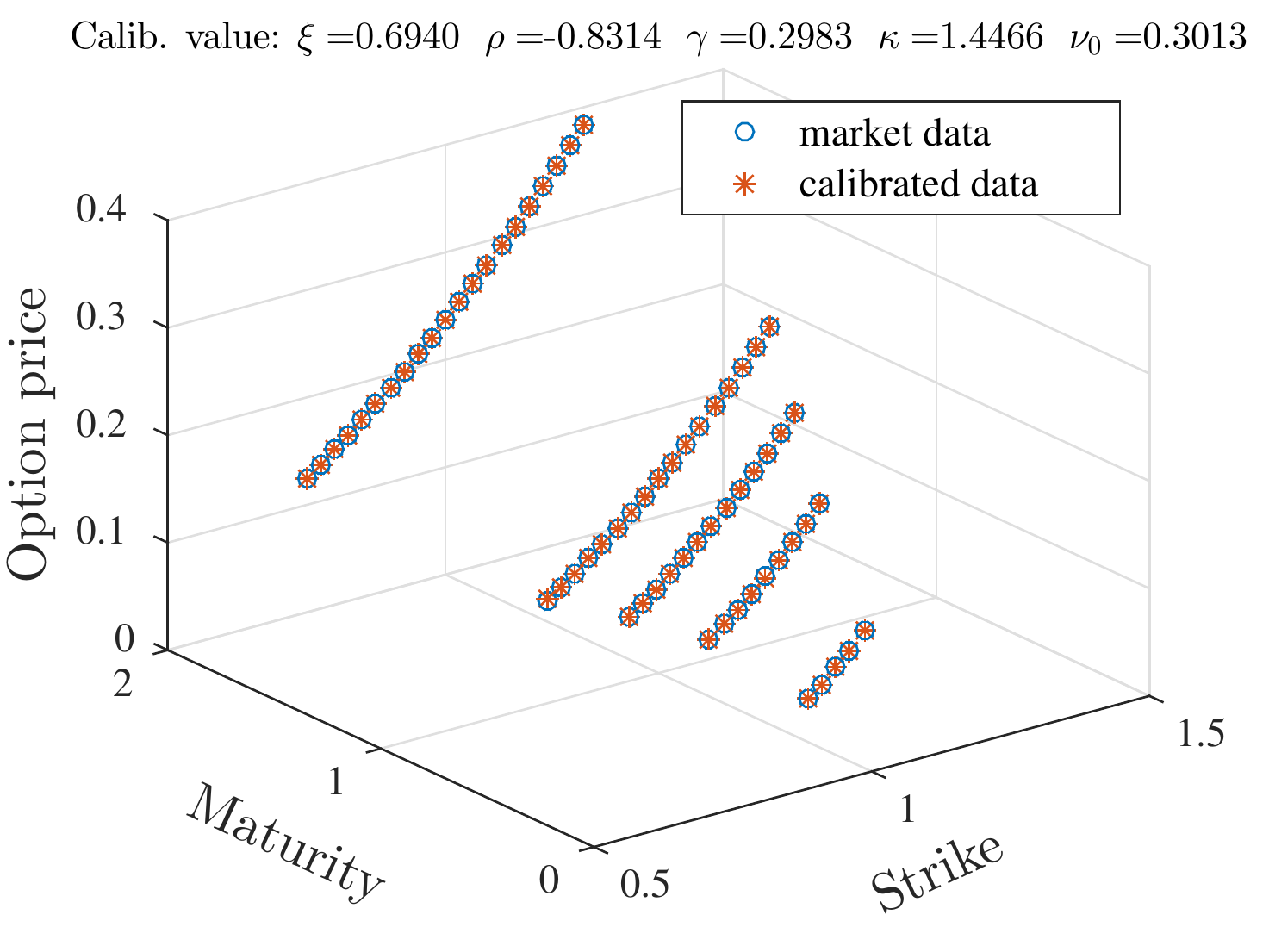} \\
 \includegraphics[width=.45\linewidth]
 {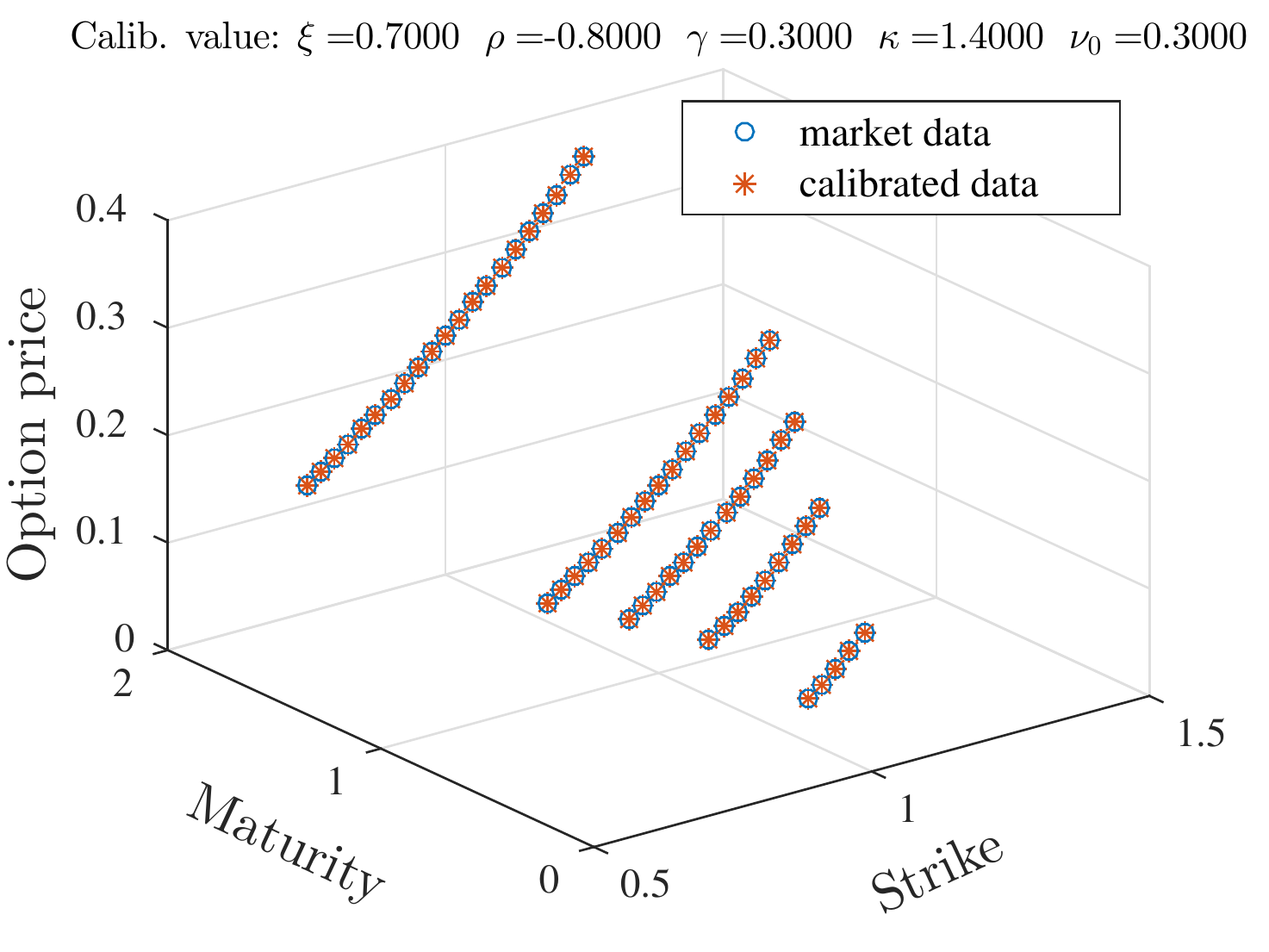}
 \includegraphics[width=.45\linewidth]
 {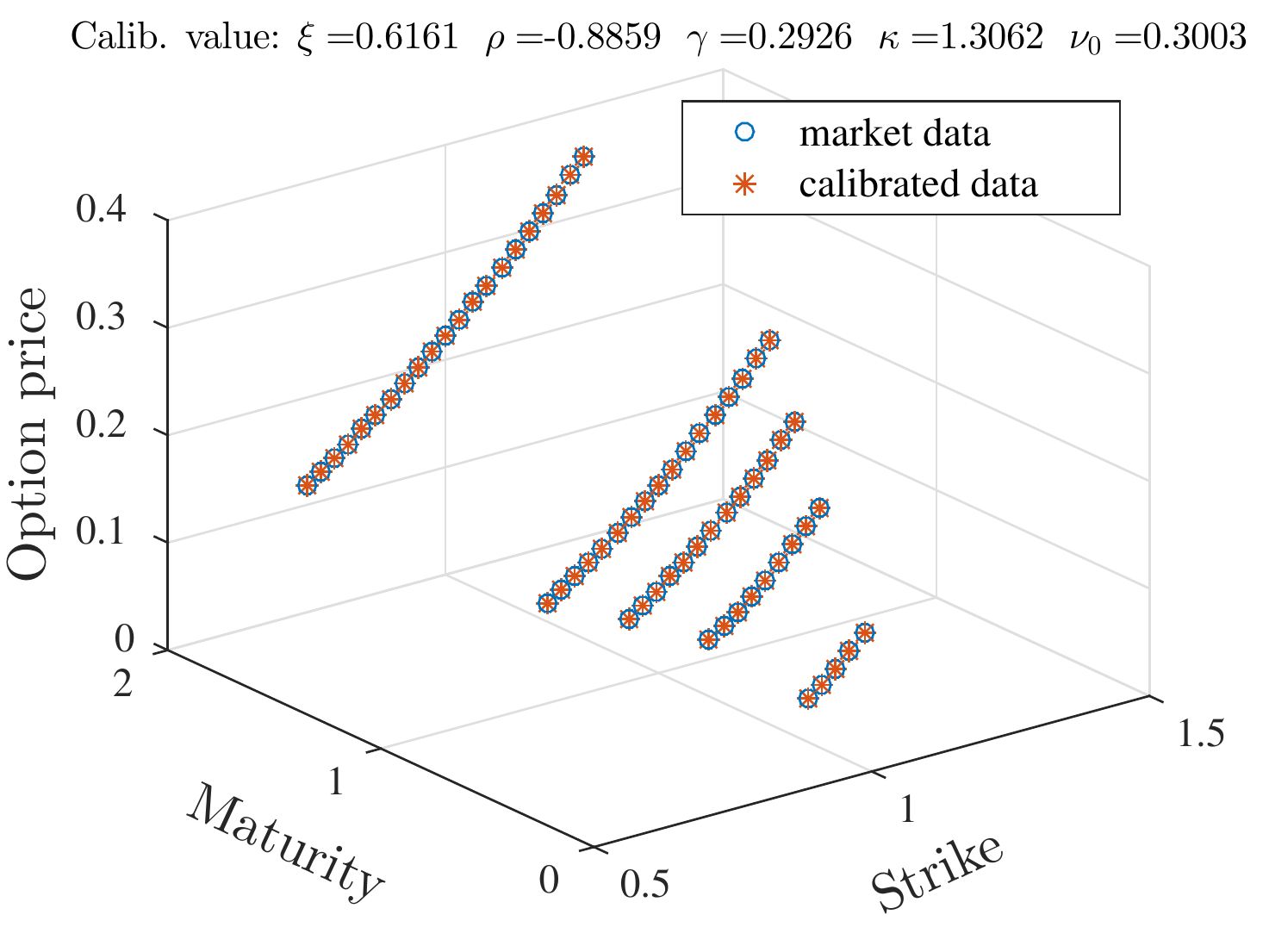}                
 \caption{Results of the calibration to the synthetic data set of American (upper row) and European (lower row) put options in the Heston model using the detailed model $\mathbb{E}_\N(\parmu)$ (left)
 and the reduced surrogate model $\mathbb{E}_N(\parmu)$ (right). The circles are the synthetic prices and the stars are the prices in the calibrated model. } 
\label{ch6:Fig:calib_ao_syn}
\end{figure}

The influence of the calibration process on the accuracy of the option price for different strike and maturity values  are shown in Figure~\ref{ch6:Fig:calib_ao_syn}  
for both American and European options.
In all plots, we hardly observe any differences in the option price obtained from the synthetic and the calibrated data with the detailed problem and the reduced problem.

 %\todo{legend of the figures}
The run-time performance is reported in Table~\ref{ch6:calin_eo_ao_syn}.
Additionally, we state the number of required iteration steps in the optimization algorithm and the number of function calls.

\begin{table}[ht]
%\begin{subtable}[c]{\textwidth}
\centering
\begin{tabular}{cccccc}                         
\toprule                                                 
 Method& $\mathbb{E}(\parmu)$&  \# iter.& \# $J$ & calib. time&
$J(\muopt^{\star})$ \\    
\midrule                                       
 $J_\N(\muopt)$&$\mathbb{E}_\N^{\rm Am}$  &7 &48 &15.59 hrs&1.698e-16 \\
 $J_N(\muopt)$ &$\mathbb{E}_N^{\rm Am}$& 8 &54 &9.50 min&9.515e-9\\
 $J_\N(\muopt)$&$\mathbb{E}_\N^{\rm Eu}$ &7 &48 &11.67 hrs &1.242e-16 \\
$J_N(\muopt)$&$\mathbb{E}_N^{\rm Eu}$ &11 & 72& 1.996 min&1.157e-08 \\
\bottomrule\end{tabular}
\caption{Calibration results for the synthetic data set of American and European put options in the Heston model in terms of the run-time performance. The number ``\# iter.'' corresponds to the number of iterations and ``\# $J$'' is the total number of function evaluations performed by the optimization routine.}
\label{ch6:calin_eo_ao_syn}
\end{table}

In this example, the optimization routine with the surrogate reduced model is about $100$
times faster for American put options and about $350$ times faster for European
put options.
The reduced model for American options
recovers the parameter slightly better than the European one. This fact can be explained by the
larger dimension of the reduced system for American options and the larger
training set, which is also reflected in the run-time performance. We point out,
that
depending on the priority of the task, i.e., accuracy vs. the run-time, one can
always manually adjust the dimension of the reduced system.

\begin{figure}[ht]
\centering
 \includegraphics[width=.45\linewidth]
 {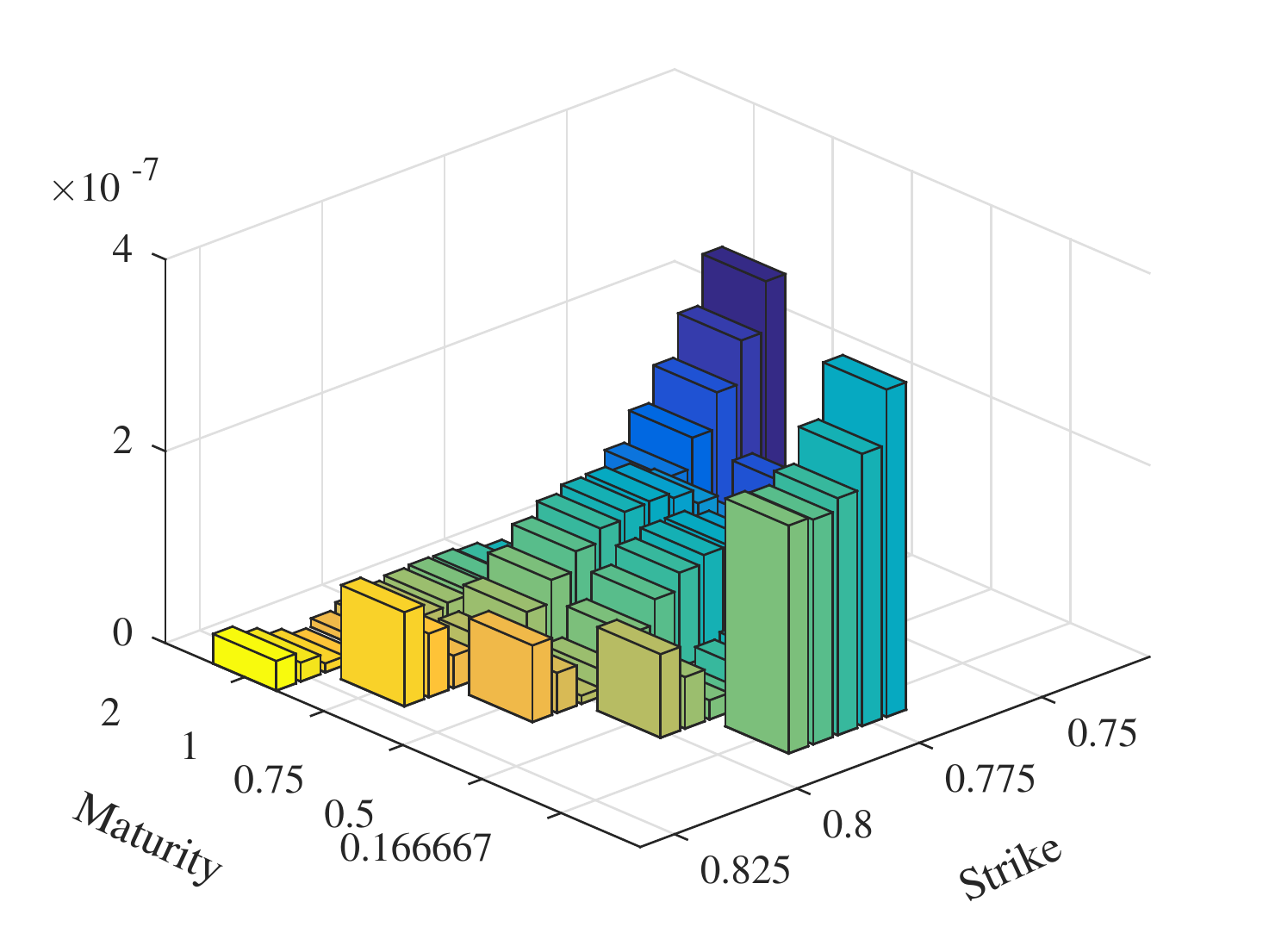}
 \includegraphics[width=.45\linewidth]
                 {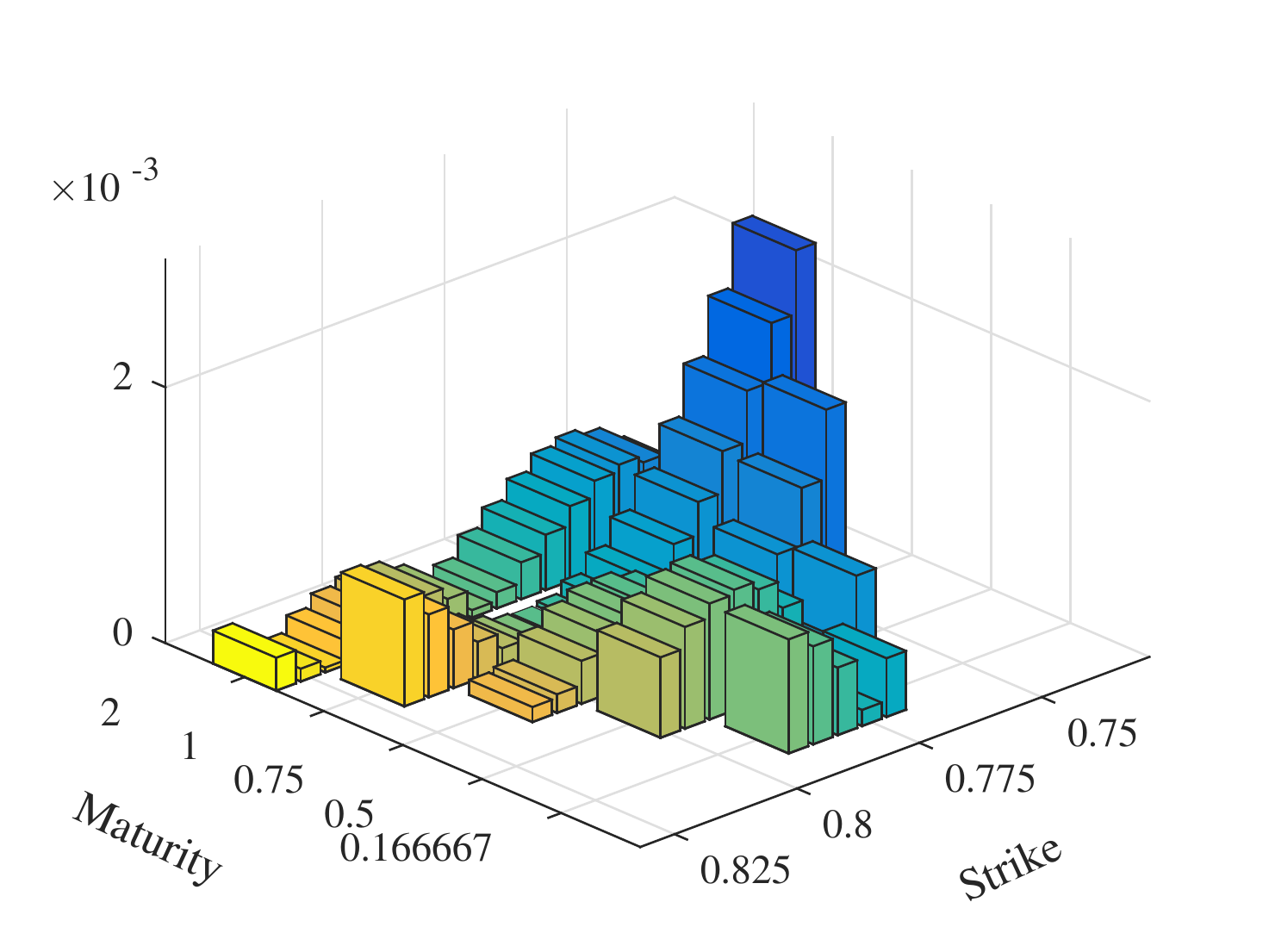}\\
                 \includegraphics[width=.45\linewidth]
 {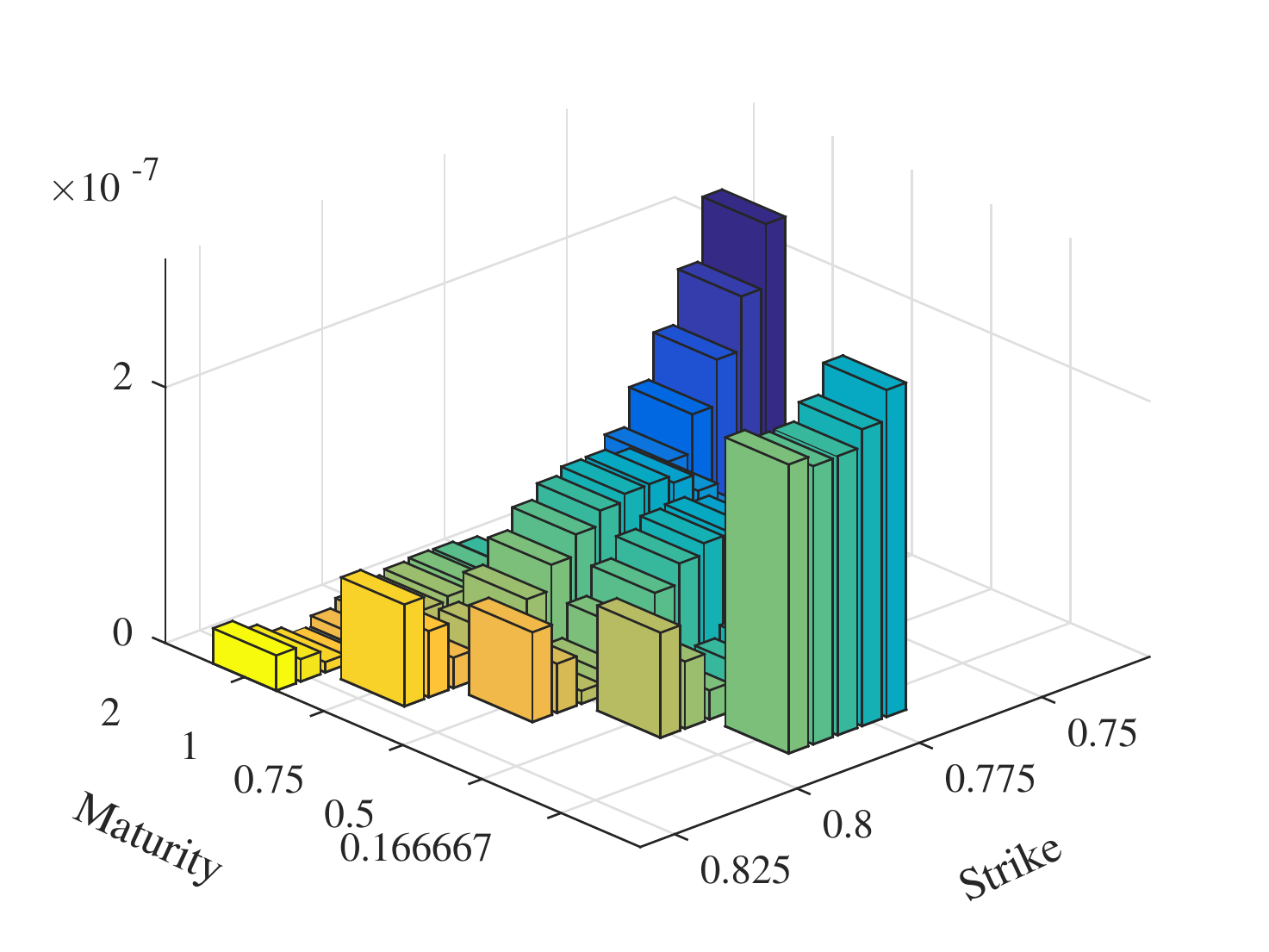}
 \includegraphics[width=.45\linewidth]
 {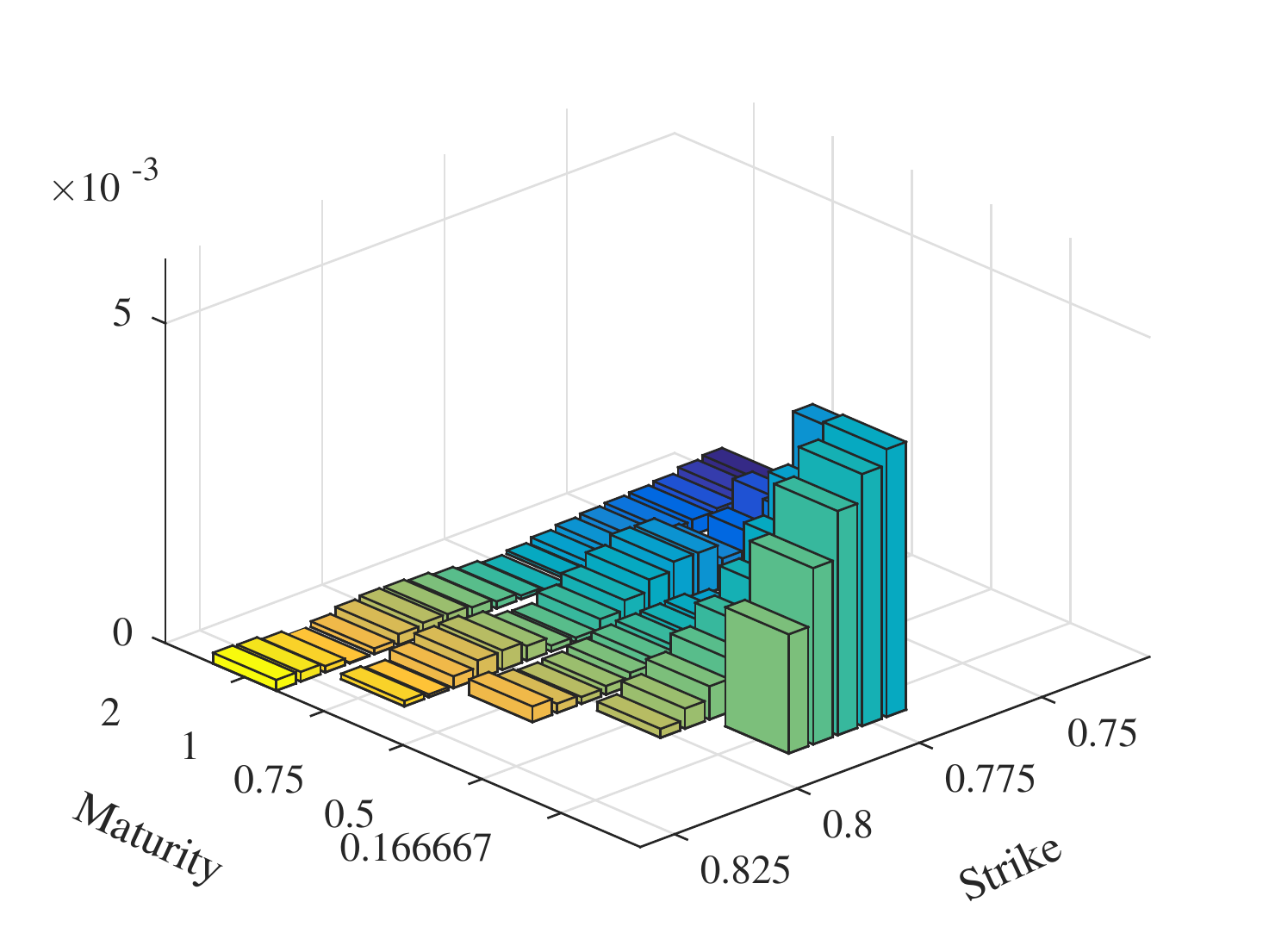}
 \caption{Calibration results for the synthetic data set of American (upper row) and European (lower row)  put options in terms of point-wise
 absolute relative errors, $|P_i^{\rm obs}-P_i^{s,\rm Am}(\muoptstar)|/P_i^{\rm obs}$, $i=1,\dots,M$, using the Heston model. 
 Left: calibration with the detailed model 
 $\mathbb{E}_\N(\parmu)$, $s=\N$. Right: calibration with the reduced surrogate model $\mathbb{E}_N(\parmu)$, $s=N$.} 
\label{ch6:Fig:calib_ao_syn_mu}
\end{figure}

To quantify the differences between reduced and detailed calibration, we plot the point-wise relative
errors $|P_i^{\rm obs}-P_i^{s}(\muoptstar)|/P_i^{\rm obs}$, $s=\{\N,N\}$, $i=1,\dots,M$ 
in Figure~\ref{ch6:Fig:calib_ao_syn_mu}.
We observe that the reduced models yield a very good fit to the synthetic data with  relative errors within the $0.5\%$-margin.

\subsection{Influence of the DAS on option pricing}\label{ch6:sec:num:deam} 
To study the effects of the DAS, we consider the difference between
the de-Americanized put option prices and the corresponding European ones. We set $S_0=1$, $r=5\%$,
\begin{align*}
K&=\left\{0.80,0.85,0.90,0.95,1.00,1.05,1.10,1.15,1.20\right\}\notag\\
T&=\left\{\frac{1}{12},\frac{2}{12},\frac{3}{12},\frac{4}{12},\frac{6}{12},\frac{9}{12}
,\frac{12}{12},\frac{24}{12}\right\}
\label{DeAm_Pricing_Setting}
\end{align*}
and consider five different  parameter sets, presented in Table~\ref{ch6:param_deam}.
 The leverage effect, that is, negative correlation between stocks and their volatilities, is commonly observed in statistical tests.
This effect was cited as one of the motivating reasons for introducing local volatility models, see \cite{cox1975notes}. In our study, we therefore consider only negative values of the correlation parameter $\rho$. We note that in some economical studies, e.g., \cite{figlewski2000leverage}, it was shown that the negative correlations between stocks and their volatilities sometimes are not primarily connected to the actual leverage, i.e., capital structure of a company, but to the general behavior in periods when the market is in downward movement.

 \begin{table}[ht]
\centering
\begin{tabular}{cccccc}                       
\toprule                                          
 Scenario&$\xi$&$\rho$ &$\gamma$ &$\kappa$ &$\vol_0$\\
 \midrule
 $p_1$&0.10&-0.20&0.07&0.1&0.07\\
 $p_2$&0.25&-0.50&0.10&0.4&0.10\\
 $p_3$&0.40&-0.50&0.15&0.6&0.15\\
 $p_4$&0.55&-0.45&0.20&1.2&0.20\\
 $p_5$&0.70&-0.80&0.30&1.4&0.30\\
\bottomrule                                         
\end{tabular} 
\caption{Overview of the parameter sets to study the DAS effect}
\label{ch6:param_deam}
\end{table}

Starting from scenario $p_1$, we increase both the
volatility of volatility parameter $\xi$ and the correlation $\rho$, and
accordingly increase the long-run  variance $\gamma$ and the mean reversion
speed $\kappa$. In all scenarios, the initial volatility is set to the value of the
mean reversion $\gamma$, i.e., $\vol_0=\gamma$. 
Using these settings, we simulate the respective American put option prices
in the Heston model, $P_i^{\N,\rm Am}$, $i=1,\dots,72$. Then, applying the
de-Americanized approach, we
translate these option prices into the corresponding European put options
$\deam{P}_i^{\N,\rm Eu}:=\tdeam(P_i^{\N,\rm Am})$. 
To
investigate the error produced by this transformation, we compare the European
put options produced by the de-Americanization method $\deam{P}_i^{\N,\rm Eu}$
with the European put
options obtained by solving the Heston PDE directly, $P_i^{\N,\rm Eu}$. 
The
values of the maximal error 
($\max_i|\deam{P}_i^{\N,\rm Eu}(\muopt)-P_i^{\N,\rm Eu}(\muopt)|$) are presented
in Table~\ref{ch6:Tab_Pricing_Heston}. The reference values
($\max_{i}P_i^{\N,\rm Eu}(\muopt)$) are presented in the third column.

\begin{table}[ht]
\centering
\setlength{\tabcolsep}{1mm}\small{
\begin{tabular}{c|rrrrrrrr|rrrrrrrr}
\toprule
&\multicolumn{8}{c|}{Maximal absolute difference}&\multicolumn{8}{c}{Maximal European price}\\
&$T_1$&$T_2$&$T_3$&$T_4$&$T_5$&$T_6$&$T_7$&$T_8$&$T_1$&$T_2$&$T_3$&$T_4$&$T_5$&$T_6$&$T_7$&$T_8$\\
	\midrule
\multirow{1}{*}{$p_1$}&		5e-4	&	2e-3	&	2e-4	&	1e-4	&	2e-4	&	2e-4	&	9e-5	&	3e-4	&	0.195	&	0.193	&	0.191	&	0.191	&	0.192	&	0.194	&	0.195	&	0.198	\\
\multirow{1}{*}{$p_2$}
&		8e-4	&	6e-5	&	3e-4	&	9e-5	&	2e-5	&	2e-5	&	7e-5	&	4e-4	&	0.196	&	0.195	&	0.195	&	0.197	&	0.200	&	0.204	&	0.208	&	0.217	\\	
\multirow{1}{*}{$p_3$}
&		2e-4	&	1e-4	&	1e-4	&	1e-5	&	1e-4	&	3e-4	&	4e-4	&	1e-3	&	0.197	&	0.199	&	0.203	&	0.207	&	0.215	&	0.224	&	0.231	&	0.250	\\	
\multirow{1}{*}{$p_4$}
&		3e-4	&	6e-5	&	7e-5	&	8e-5	&	2e-4	&	3e-4	&	4e-4	&	1e-3	&	0.199	&	0.205	&	0.211	&	0.218	&	0.229	&	0.243	&	0.255	&	0.286	\\
\multirow{1}{*}{$p_5$}	
&		2e-4	&	1e-4	&	1e-4	&	3e-5	&	1e-4	&	4e-4	&	9e-4	&	3e-3	&	0.202	&	0.213	&	0.223	&	0.233	&	0.249	&	0.268	&	0.283	&	0.326	\\
\bottomrule
\end{tabular}}
\caption{Effects of the DAS on put option pricing in the Heston model.}
\label{ch6:Tab_Pricing_Heston}
\end{table}

We observe especially for scenarios $p_4$ and $p_5$, where the volatility of
volatility $\xi$ and
correlation $\rho$ have higher
values, that the DAS has the strongest effect. {Note that here we fix $r=5\%$ as a particular case to highlight the de-Americanization effect, and we refer to \cite{deam2016} for a more detailed parameter study.}

\subsection{Comparison between RBM and DAS}\label{ch6:sec:comparison}
Next, we perform a numerical comparison of the calibration with
American put options using both model reduction
techniques.
That is, we consider the minimization
problems~\eqref{ch6:minprob_rb} and~\eqref{ch6:minprob_deam}.
For comparative purposes, we also carry out
the calibration with the detailed finite element solution \eqref{ch6:minprob}.

First, we use a synthetic set of observations $P^{\rm obs}$  given by~\eqref{eq:syndatasetting}. We consider different parameter scenarios
corresponding to different values of $\muopt\in\popt$; see
Table~\ref{ch6:Table_data_set}. For each scenario, we construct an artificial
set of observations $P_i^{\rm obs}:=P_i^{\N,\rm Am}$, $i=1,\dots,65$. In
general, the parameter $\kappa$ is
price-insensitive, see, e.g., \cite{janek2011fx}, and thus it cannot be reconstructed.  Therefore, we fix $\kappa$ to its exact value and
do the parameter estimation only for $\xi,\rho,\gamma$, and $\vol_0$.  
The results of the calibration are summarized in
Tables~\ref{ch6:calib_ao_syn_data} and
\ref{ch6:calib_ao_syn_data_1}. We observe that, overall, all methods provide
a good reconstruction of the parameters. As can be expected, the most cost-intense variant also provides the results with the highest accuracy.
However all of our surrogate models yield quite good results.
\begin{table}[ht!]
\centering
\begin{tabular}{ccccccc}                       
\toprule                                     
 &$\xi$&$\rho$ &$\gamma$ &$\kappa$ &$\vol_0$\\
 \midrule
 $p_1$&0.10&-0.20&0.07&0.5&0.07\\
 $p_2$&0.25&-0.50&0.10&0.5&0.10\\
 $p_3$&0.40&-0.50&0.15&0.6&0.15\\
 $p_4$&0.55&-0.45&0.20&1.2&0.20\\
 $p_5$&0.70&-0.80&0.30&1.4&0.30\\
 $p_6$&0.2928&-0.7571&0.0707&0.6067&0.0707\\ 
\bottomrule                                            
\end{tabular} 
\caption{Overview of the parameter sets used to generate the
synthetic data set}
\label{ch6:Table_data_set}
\end{table}

\begin{table}[ht]
\centering\begin{small}
\begin{tabular}{clccccc}                         
\toprule                                                 
 Scenario&Method &$\mathbb{E}(\parmu)$ &  $\# J$ & calib. time& pre-process. time for $P^{\rm obs}$\\
\midrule                                           
&$J_\N(\muopt)$  &$\mathbb{E}_\N^{\rm Am}$ &75 &8.524 hrs &\\
$p_1$& $\widetilde{J}_\N(\muopt)$  &$\mathbb{E}_\N^{\rm Eu}$ &145 & 18.018 min&36.045 min\\
&${J}_N(\muopt)$  &$\mathbb{E}_N^{\rm Am}$ & 70& 3.912 min& \\
\midrule
%& RB DeAm&170&  4.855 min&36.045 min\\
&${J}_\N(\muopt)$  &$\mathbb{E}_\N^{\rm Am}$ & 65& 8.035 hrs&\\
$p_2$&$\widetilde{J}_\N(\muopt)$  &$\mathbb{E}_\N^{\rm Eu}$ &70 & 8.895 min&36.830 min\\
&${J}_N(\muopt)$  &$\mathbb{E}_N^{\rm Am}$ & 70& 3.489 min&  \\
\midrule
%& RB DeAm&70& 2.050 min&36.830 min \\
&${J}_\N(\muopt)$  &$\mathbb{E}_\N^{\rm Am}$ & 60&8.341 hrs &\\
$p_3$&$\widetilde{J}_\N(\muopt)$  &$\mathbb{E}_\N^{\rm Eu}$ &60 &   8.083 min&35.374 min \\
&${J}_N(\muopt)$  &$\mathbb{E}_N^{\rm Am}$ &70 & 3.574 min&  \\
%& RB DeAm&80& 2.353 min& 35.374 min\\
\midrule
&${J}_\N(\muopt)$  &$\mathbb{E}_\N^{\rm Am}$ & 50&7.088 hrs &\\
$p_4$&$\widetilde{J}_\N(\muopt)$  &$\mathbb{E}_\N^{\rm Eu}$ &45 &  6.031 min&36.660 min \\
&${J}_N(\muopt)$  &$\mathbb{E}_N^{\rm Am}$ &55 & 2.813 min&  \\
%& RB DeAm&75&  2.189 min& 36.660 min\\
\midrule
&${J}_\N(\muopt)$  &$\mathbb{E}_\N^{\rm Am}$  &40 &6.331 hrs &\\
$p_5$&$\widetilde{J}_\N(\muopt)$  &$\mathbb{E}_\N^{\rm Eu}$ &65 & 8.607 min&36.574 min \\
&${J}_N(\muopt)$  &$\mathbb{E}_N^{\rm Am}$ &45 & 2.271 min&  \\
%& RB DeAm&65&  1.933 min&36.574 min \\
\midrule
&${J}_\N(\muopt)$  &$\mathbb{E}_\N^{\rm Am}$  & 70& 9.665 hrs&\\
$p_6$&$\widetilde{J}_\N(\muopt)$  &$\mathbb{E}_\N^{\rm Eu}$ & 70& 9.555 min&36.960 min \\
&${J}_N(\muopt)$  &$\mathbb{E}_N^{\rm Am}$ &90 & 4.503 min&  \\
%& RB DeAm&100&3.055 min & 36.960 min\\
\bottomrule
\end{tabular}   
\caption{Computational time for calibrating American put options
  using different model reduction techniques.}
   \label{ch6:calib_ao_syn_data}  \end{small}                                           
\end{table}

We observe that both reduction approaches provide a significant speed-up compared
to the expensive detailed solver, which on average takes about eight hours for each
scenario. Although the DAS allows for an extremely 
efficient calibration process, it requires an additional pre-processing step for
the data. In contrast to the RBM, this pre-processing depends on the actual market data and therefore can not be performed in advance. This
is a serious bottleneck compared to the RBM approach.

\begin{table}[ht!]
\centering
\begin{small}
\begin{tabular}{cllllllll}                         
\toprule                                               
Scenario&Method &$\mathbb{E}(\parmu)$&$\muopt$&$\xi$ &$\rho$ &$\gamma$ &$\vol_0$ &$J(\muopt^\star)$\\
\midrule
&&&$\muopt_{\rm ex}$&0.1&-0.2 &0.07  &0.07 &\\
$p_1$&${J}_\N(\muopt)$  &$\mathbb{E}_\N^{\rm Am}$& $\muopt^{\star}$&0.1002 &-0.1997 &0.07  &0.07  &7.8806e-14\\
&$\widetilde{J}_\N(\muopt)$  &$\mathbb{E}_\N^{\rm Eu}$&$\muopt^{\star}$&0.1000   &-0.3003    &0.0701    &0.0688& 1.0248e-07\\
&${J}_N(\muopt)$  &$\mathbb{E}_N^{\rm Am}$&$\muopt^{\star}$&0.1477  & -0.0788&    0.0697&    0.0700  &4.3440e-08\\
%&RB DeAm&0.1849 &  -0.1801   & 0.0772   & 0.0634&3.2571e-07\\
\midrule
&&&$\muopt_{\rm ex}$&0.25 &-0.5 &0.1 &0.1 &\\
$p_2$&${J}_\N(\muopt)$  &$\mathbb{E}_\N^{\rm Am}$&$\muopt^{\star}$&0.25 &-0.5 &0.1  &0.1  &6.4756e-17\\
&$\widetilde{J}_\N(\muopt)$  &$\mathbb{E}_\N^{\rm Eu}$&$\muopt^{\star}$& 0.2404   &-0.5388    &0.1001   & 0.0991&1.1363e-08\\
&${J}_N(\muopt)$  &$\mathbb{E}_N^{\rm Am}$&$\muopt^{\star}$&0.2860   &-0.4824    &0.0968   & 0.1015&9.9148e-08\\
%&RB DeAm&0.2973  & -0.4104    &0.1087    &0.0959&1.1681e-07\\
\midrule
&&&$\muopt_{\rm ex}$&0.4 &-0.5&0.15  &0.15 &\\
$p_3$&${J}_\N(\muopt)$  &$\mathbb{E}_\N^{\rm Am}$&$\muopt^{\star}$&0.4 &-0.5 &0.15  &0.15 &2.9834e-18\\
&$\widetilde{J}_\N(\muopt)$  &$\mathbb{E}_\N^{\rm Eu}$& $\muopt^{\star}$&0.4282 &  -0.4620   & 0.1544   & 0.1492& 2.2003e-09\\
&${J}_N(\muopt)$  &$\mathbb{E}_N^{\rm Am}$& $\muopt^{\star}$&0.3537   &-0.5731    &0.1456    &0.1504&6.2085e-09\\
%&RB DeAm&  0.4576&   -0.4155   & 0.1519  &  0.1484&4.4351e-08\\
\midrule
&&&$\muopt_{\rm ex}$&0.55 &-0.45 &0.2 &0.2 &\\
$p_4$&${J}_\N(\muopt)$  &$\mathbb{E}_\N^{\rm Am}$&$\muopt^{\star}$&0.5502 &-0.4499 &0.2   &0.2  &4.8235e-14  \\
&$\widetilde{J}_\N(\muopt)$  &$\mathbb{E}_\N^{\rm Eu}$&$\muopt^{\star}$& 0.5801  & -0.4220  &  0.2044 &   0.1989&1.5377e-09\\
&${J}_N(\muopt)$  &$\mathbb{E}_N^{\rm Am}$& $\muopt^{\star}$&0.5048   &-0.4980   & 0.1989  &  0.1995& 1.6681e-08\\
%&RB DeAm&   0.5359   &-0.4473  &  0.2013 &   0.1989&1.9256e-08\\
\midrule
&&&$\muopt_{\rm ex}$&0.7 &-0.8 &0.3 &0.3 &\\
$p_5$&${J}_\N(\muopt)$  &$\mathbb{E}_\N^{\rm Am}$& $\muopt^{\star}$&0.7 &-0.8 &0.3   &0.3  &3.4388e-18  \\
&$\widetilde{J}_\N(\muopt)$  &$\mathbb{E}_\N^{\rm Eu}$&$\muopt^{\star}$&  0.8433 &  -0.6668  &  0.3170 &   0.2990& 1.8369e-08\\
&${J}_N(\muopt)$  &$\mathbb{E}_N^{\rm Am}$& $\muopt^{\star}$&0.6881&   -0.8259  &  0.2994  &  0.3006& 1.0533e-08\\
%&RB DeAm&    0.7718 &  -0.7136 &   0.3102 &   0.2991&2.3145e-08\\
\midrule
&&&$\muopt_{\rm ex}$&0.2928 &-0.7571 &0.0707 &0.0707 &\\
$p_6$&${J}_\N(\muopt)$  &$\mathbb{E}_\N^{\rm Am}$&$\muopt^{\star}$&0.2928  &0.7571 &0.0707  &0.0707  &6.5746e-18  \\
&$\widetilde{J}_\N(\muopt)$  &$\mathbb{E}_\N^{\rm Eu}$&$\muopt^{\star}$&  0.3690   &-0.6026   & 0.0736  &  0.0685&1.6179e-07\\
&${J}_N(\muopt)$  &$\mathbb{E}_N^{\rm Am}$&$\muopt^{\star}$& 0.3096   &-0.7049   & 0.0700  &  0.0718&9.6369e-08 \\
%&RB DeAm& 0.3527  & -0.5209  &  0.0814  &  0.0638&3.5672e-07\\
\midrule
\end{tabular}  \end{small}
\caption{Calibration results on the synthetic data set of American put options in the Heston model with
different model reduction techniques: the RBM, the de-Americanization method and the detailed problem.}
   \label{ch6:calib_ao_syn_data_1}                                           
\end{table}

\begin{figure}[ht]
\centerline{\begin{minipage}[b]{0.45\textwidth}
 \includegraphics[width=\linewidth]{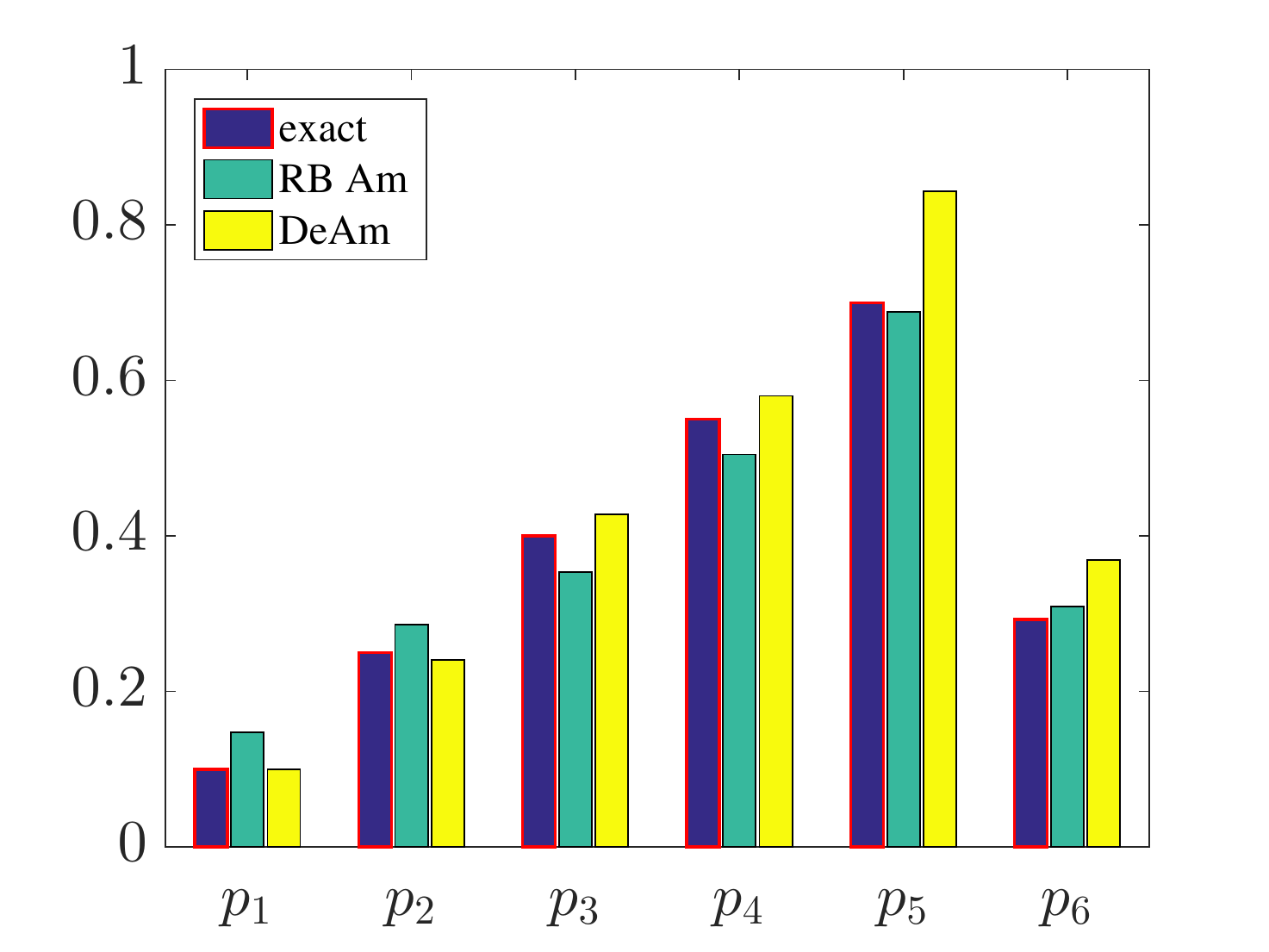} 
 \centering{$\xi$}
\end{minipage}
\begin{minipage}[b]{0.45\textwidth}
 \includegraphics[width=\linewidth]{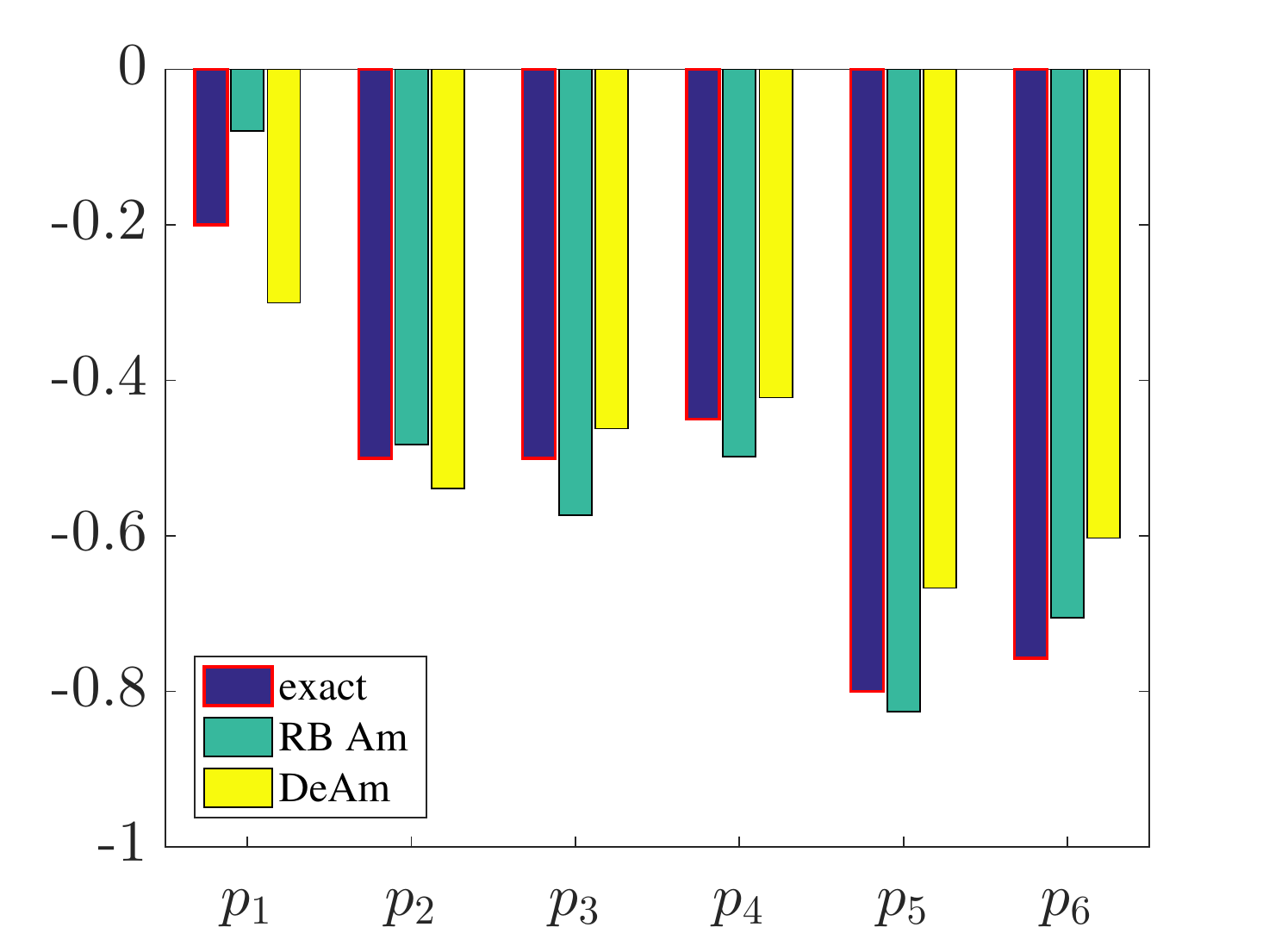} 
 \centering{$\rho$}
\end{minipage}}
\centerline{
\begin{minipage}[b]{0.45\textwidth}
 \includegraphics[width=\linewidth]{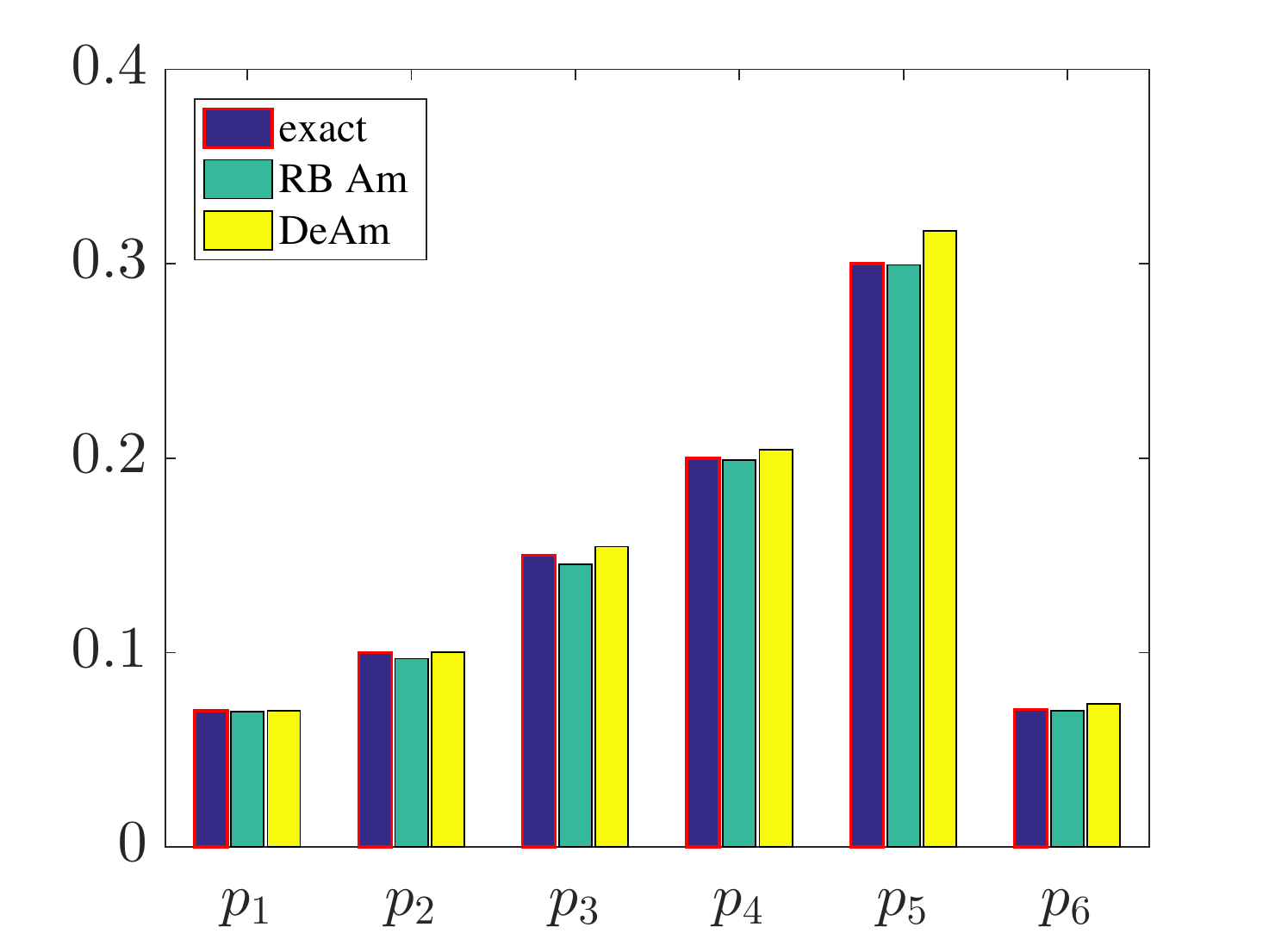} 
 \centering{$\gamma$}
\end{minipage}
\begin{minipage}[b]{0.45\textwidth}
 \includegraphics[width=\linewidth]{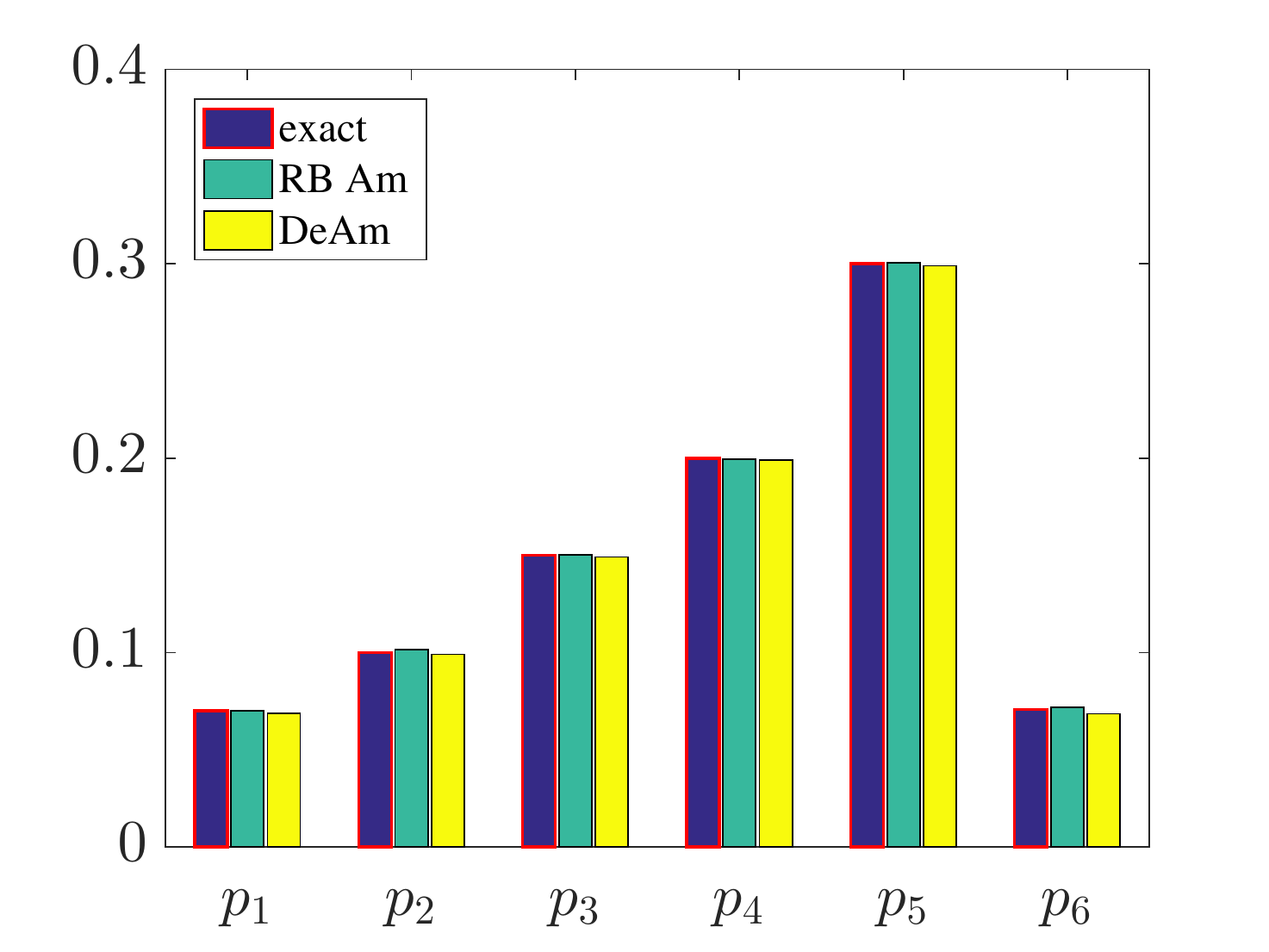} 
 \centering{$\vol_0$}
\end{minipage}}
 \caption{Reconstructed parameters for the different
scenarios obtained by
calibrating  American put options  with
different model reduction techniques.} 
\label{ch6:Fig:calib_mu}
\end{figure}

Figure~\ref{ch6:Fig:calib_mu} shows the influence of the different calibration approaches on the parameters.
It can be seen that in all approaches
the main difficulty arises in in identifying $\xi$ and $\rho$. In fact, this tendency has been also observed for the detailed solver 
(see cases $p_1$, $p_4$, Table~\ref{ch6:calib_ao_syn_data_1}). The
remaining parameters $\gamma$ and $\nu_0$ are recovered almost exactly. We also
note that for scenarios $p_1$--$p_3$, which correspond to the cases when
$\xi,\vol_0$  and $\rho$ are the smallest, the DAS calibration  is able to provide a better reconstruction of the
parameter $\xi$ than the RBM. By contrast, in the scenarios
$p_4$ and $p_5$, which correspond to large values
of the correlation parameter, the DAS gives poorer results for $\xi$ and $ \rho$. This is consistent with
the observation made in Section~\ref{ch6:sec:deam}, where we noticed that
the errors caused by de-Americanization are the largest in the cases where
$\xi$ and $\rho$ are large (see Table~\ref{ch6:Tab_Pricing_Heston},
scenario $p_5$, maturity $T_8$). 

To summarize our findings, the cases with ``extreme''  parameter
values have a significant impact on the performance of the optimization
routine, in both the detailed and the reduced problems. In the case of
the reduced basis method, this difficulty may be overcome for example by
considering an adaptive parameter domain partition~\cite{haasdonk2011mutrain}, in particular for $\xi$
and $\rho$, or by increasing the number of snapshots, and furthermore by increasing the dimension of the reduced system. However, this will also result in an increase in the time required to compute the solution.

\subsection{Real market data}\label{ch6:sec:google}
Finally, we extend our approach to the calibration of a real market data set, provided by options on the Google stock. 
Since the Google stock does not pay dividends, the American call options can be 
priced the same as the European call options, \cite{hull03}. 
Hence, we restrict consideration to only American put options. Namely, we consider the data $P^{\rm obs}$ of 401 American 
put options with $S_0=523.755$, $r=0.15\%$ on February 2nd, 2015. The data is pre-processed using the methodology applied to the volatility index (VIX) by the Chicago board of 
exchange~\cite{vix}: 
\begin{itemize}
 \item For each option with strike price $K_i$, we consider the midpoint of the bid-ask spread.
 \item Options with zero bid prices are neglected.
 \item If two puts with consecutive strike prices have zero bid prices, no puts with lower strike prices are included.
\end{itemize}

The used data is given in the appendix, see  Table~\ref{ch6:table:google_data}. In
terms of the
moneyness, we consider all types, out-of-the money ($K_i<S_0$), at-the-money
($K_i=S_0$) and in-the-money ($K_i>S_0$) options. 

In our synthetic test scenarios,  the Feller condition~\eqref{ch2:feller}
was automatically satisfied.
However,  this does not hold for general calibration processes. Thus,
we impose the following additional constraint on 
$\muopt=(\xi,\rho,\gamma,\kappa,\vol_0)\in\popt$
\begin{equation}
 \popt:=\left\{\muopt\in\mathbb{R}^5:
\muopti_{\min,i}\leq\muopti_i\leq\muopti_{\max,i},\quad
2\muopti_3\muopti_4-\muopti_1^2<0, \quad i=1,\dots,5\right\}.
\end{equation}
As optimization algorithm, we take the MATLAB function
{\it fmincon} based on  the Interior-Point method and which, in contrast to {\it
lsqnonlin}, allows the inclusion of inequality constraints. {We consider the same termination condition for the optimization routine as previously.}

To calibrate the parameters, we consider the detailed minimization 
problem~\eqref{ch6:minprob} and, as previously, the reduced models: the RBM~\eqref{ch6:minprob_rb}, 
the DAS \eqref{ch6:minprob_deam} 
and the combination of both~\eqref{ch6:minprob_deam_rb}. For completeness, we 
also consider the calibration of the de-Americanized data using the closed-form solution of \eqref{ch6:minprob_deam_cf}.

\begin{table}[ht]
\centering
\begin{tabular}{lcc c cccc}                         
\toprule                                               
Method&$\mathbb{E}(\parmu)$& $\muopt$& $\xi$ &$\rho$ &$\gamma$ &$\kappa$& $\vol_0$\\ 
\midrule
 &&$\muopt_{\rm in}$ &0.6005 &-0.6815& 0.4867 &2.02 &0.4961 \\    
  ${J}_\N(\muopt)$  &$\mathbb{E}_\N^{\rm Am}$&$\muopt^\star$ &0.5953&   -0.7210&  0.0527&    3.3615&0.0584 \\
 ${J}_N(\muopt)$  &$\mathbb{E}_N^{\rm Am}$&$\muopt^\star$ &0.5144 &-0.7964& 0.0521&2.5906 &0.0554\\
  $\widetilde{J}_\N(\muopt)$  &$\mathbb{E}_\N^{\rm Eu}$&$\muopt^\star$ &0.4095 &-0.6818& 0.0516&1.6262 &0.0567\\
  % RB DeAm&$\muopt$ &0.3785 &-0.9469& 0.0500&4.9822 &0.0717\\
   $\widetilde{J}_{\rm CF}(\muopt)$  &&$\muopt^\star$ &0.3927 &-0.6518& 0.0580&1.4554 &0.0546\\
\hline
\end{tabular}
\caption{Parameters obtained by the calibration on American put options given on the Google stock using different methods}
\label{ch6:calib_google_mu} 
%\end{subtable}
\end{table}

The results of the
calibration are presented in Table~\ref{ch6:calib_google_mu}. As before, $\gamma$ and
$\vol_0$ can be easily identified by all our approaches.
The rate of mean
reversion $\kappa$ appears to be a non-identifiable parameter, and all models provide quite different results.
For the
remaining parameters $\xi$ and $\rho$, we observe that the DAS
 tends to underestimate the volatility of volatility $\xi$ and the
correlation $\rho$, compared to the detailed and RBM approach. 
This is clearly reflected in all models that use the perturbed de-Americanized data, i.e., $\deam{J}_\N(\muopt)$ and $\deam{J}_{\rm CF}(\muopt)$.
This is in good agreement with our observations for the synthetic data sets (see scenario $p_5$ and
$p_6$), wherefor large (absolute) values of the correlation parameters
the DAS was unable to provide a good reconstruction of the parameters 
$\xi$ and $\rho$.

\begin{table}[ht]
\centering
\begin{tabular}{lcccccc}                         
\toprule                                               
 Method&$\mathbb{E}(\parmu)$&  \# iter.& \# $J$ & calib. time&pre-process. time for $P^{\rm obs}$ \\
\midrule
 ${J}_\N(\muopt)$  &$\mathbb{E}_\N^{\rm Am}$  &35 &219&68.72 hrs& 
\\
 ${J}_N(\muopt)$  &$\mathbb{E}_N^{\rm Am}$ & 38 &260 &44.20 min& \\
 $\widetilde{J}_\N(\muopt)$  &$\mathbb{E}_\N^{\rm Eu}$ &34 &207 &56.47 min & 4.96 hrs\\
 % RB DeAm &35 &228 &37.40 min & 4.96 h.\\
 $\widetilde{J}_{\rm CF}(\muopt)$  & &43 &265 &4.30 min & 4.96 hrs\\
\bottomrule
\end{tabular}
\caption{Computational time for calibrating American put options given on the Google stock in the Heston model using different methods}
\label{ch6:calib_google_time} 
\end{table}

The results of the run-time performance of the different methods is given in 
Table \ref{ch6:calib_google_time}. 
We observe that the detailed approach is much more cost-intense than the proposed surrogate models. The cost can be drastically reduced from a couple of days to
less than an hour. Notably, we observe the substantial speed-up obtained by evaluating model prices with the
closed-form solutions. This approach appears to us the most efficient when dealing with European options.
However, taking into consideration the additional time for pre-processing the data in the DAS, 
the total time for calibration with this method
can be much slower than the calibration with American options using the RBM, depending how often the calibration has to be performed with new market data. However, the DAS pre-processing time could be sped up by implementing more advanced tree methods. We also note that, in contrast to the DAS, the RBM approach allows us to control the accuracy by increasing the dimension of the reduced spaces.

\begin{figure}[ht]
\begin{minipage}[b]{\textwidth}
  \centering
  \includegraphics[width=.45\linewidth]
{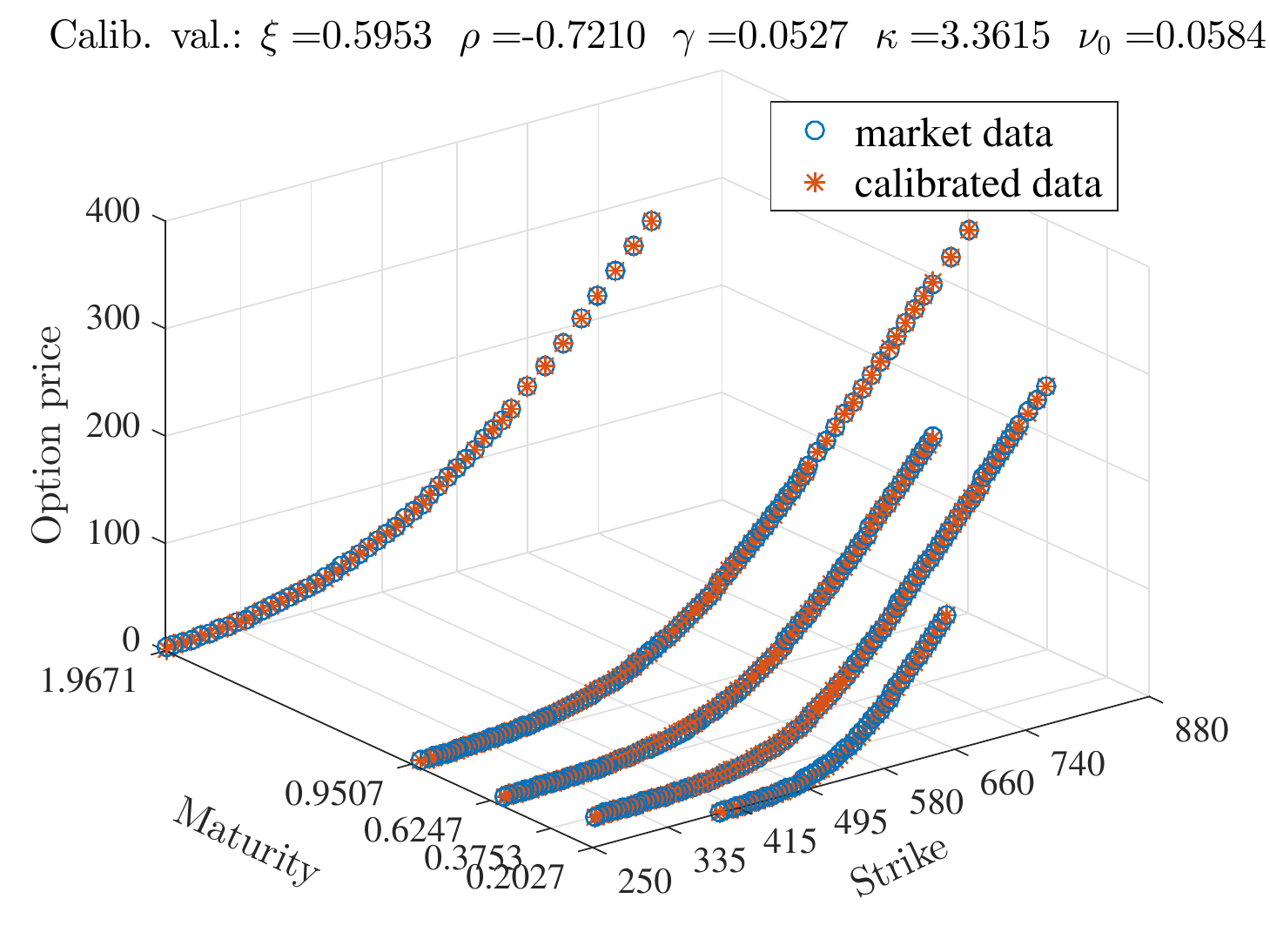}
\includegraphics[width=.45\linewidth]
{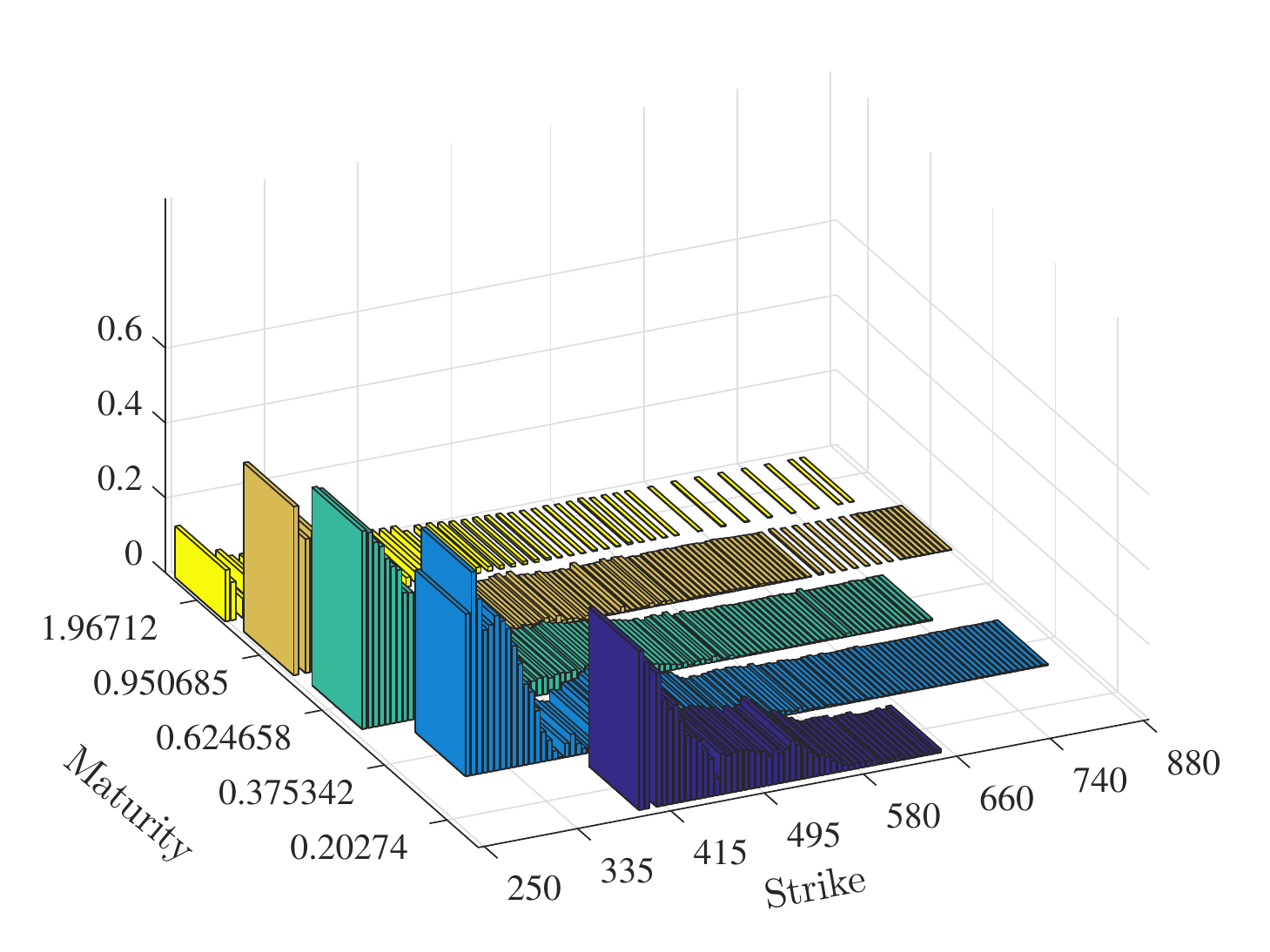}
\centering{Detailed problem, $J_\N(\muopt)$}
 \end{minipage}
 \begin{minipage}[b]{\textwidth}
   \centering
   \includegraphics[width=.45\linewidth]
{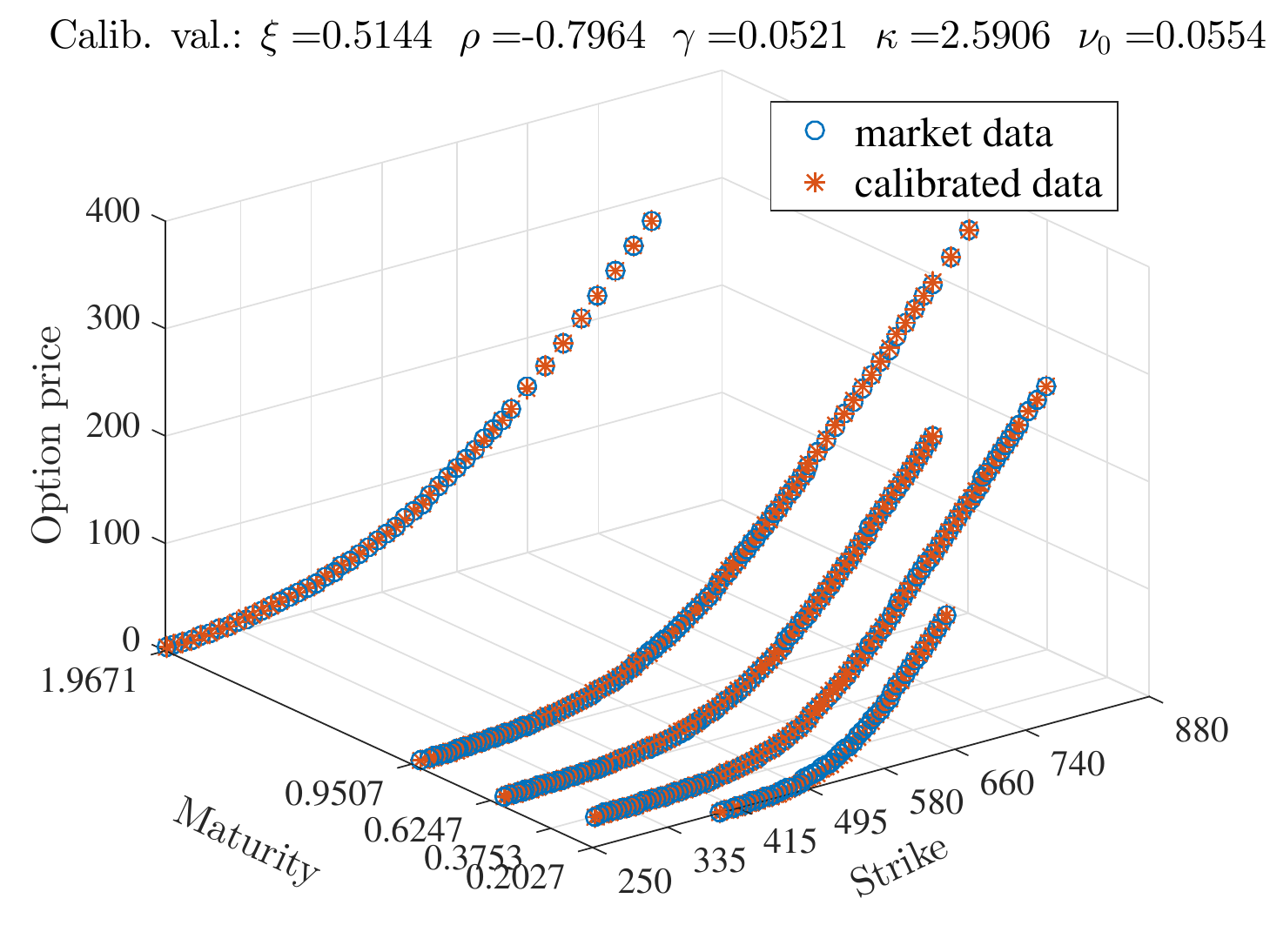}
\includegraphics[width=.45\linewidth]
{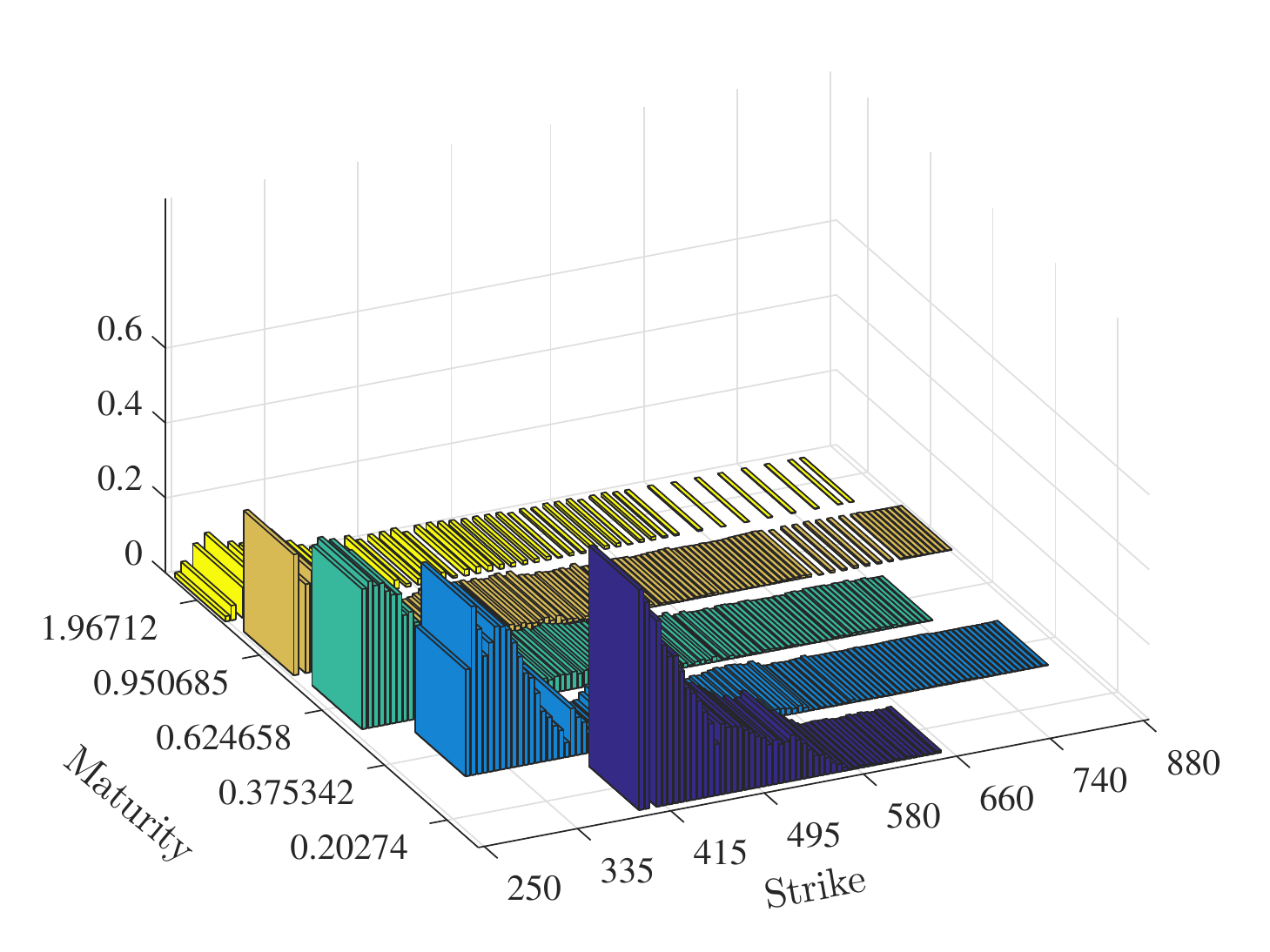}
\centering{Reduced problem, $J_N(\muopt)$}
 \end{minipage}
 \caption{Left: the Google data set of American put
options (circles) and the calibrated model data in the Heston model (stars). 
Right: the relative error of
the market and calibrated data, $|P_i^{\rm obs}-P_i^{s,\rm Am}(\muopt^\star)|/P_i^{\rm obs}$, $i=1,\dots,M$, $s=\{\N,N\}$} 
\label{ch6:fig_google_ao_actual}
\end{figure}

\begin{figure}[ht]
\centering
\begin{minipage}[b]{\textwidth}
 \includegraphics[width=.45\linewidth]{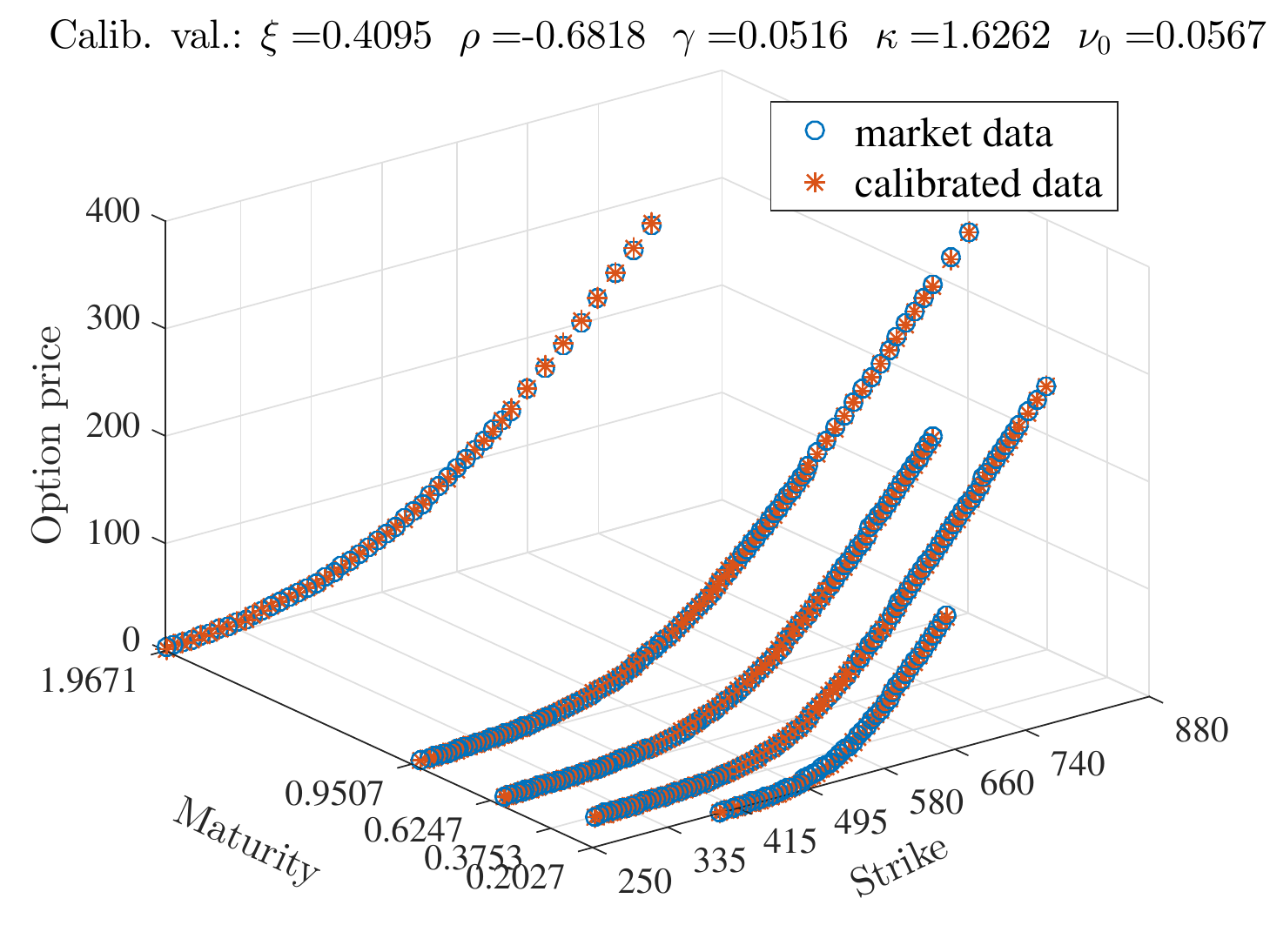}
\includegraphics[width=.45\linewidth]
{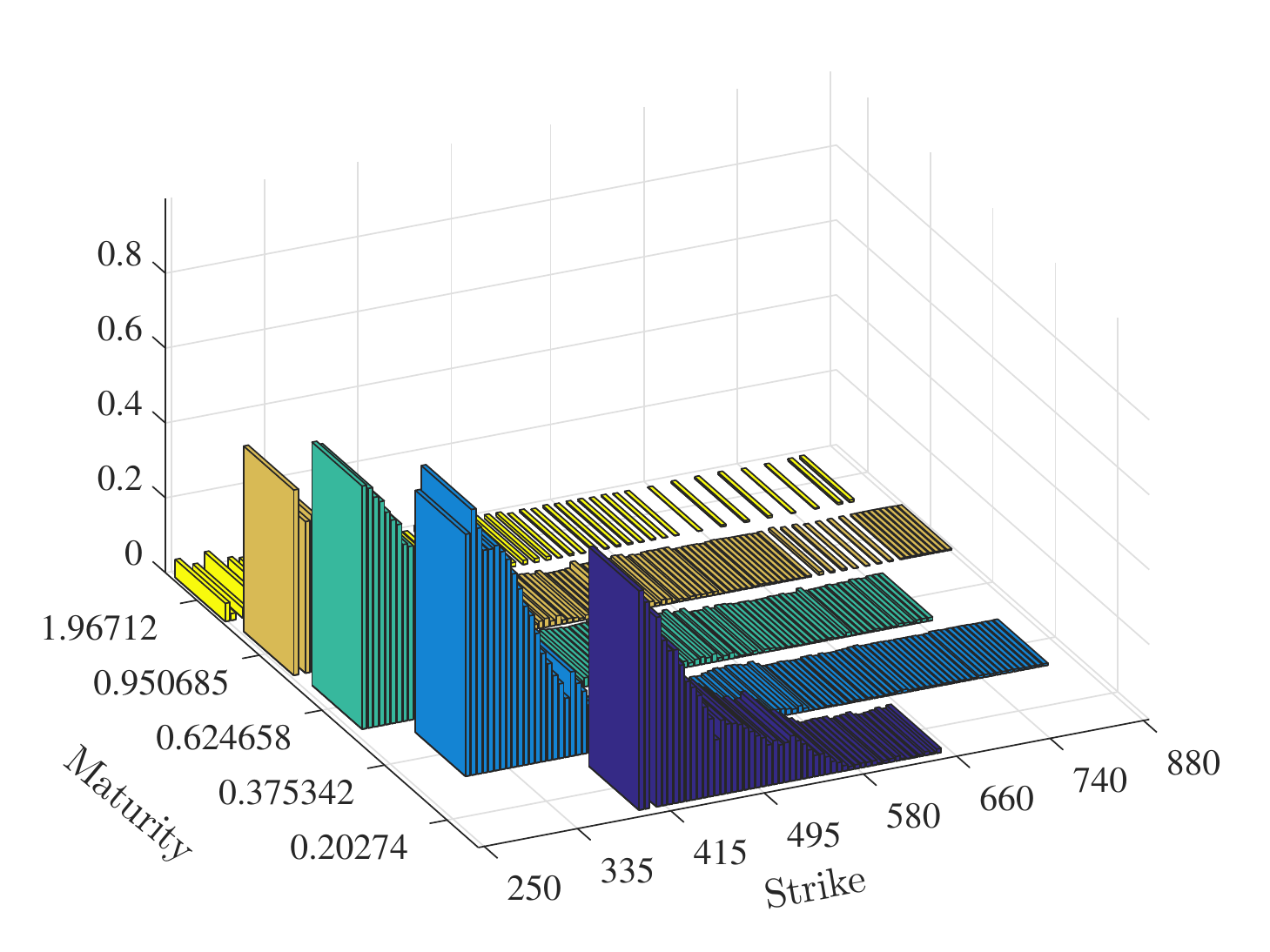}
\centering{De-Americanized problem, $\deam{J}_\N(\muopt)$}
 \end{minipage}
 
 \begin{minipage}[b]{\textwidth}
  \includegraphics[width=.45\linewidth]
{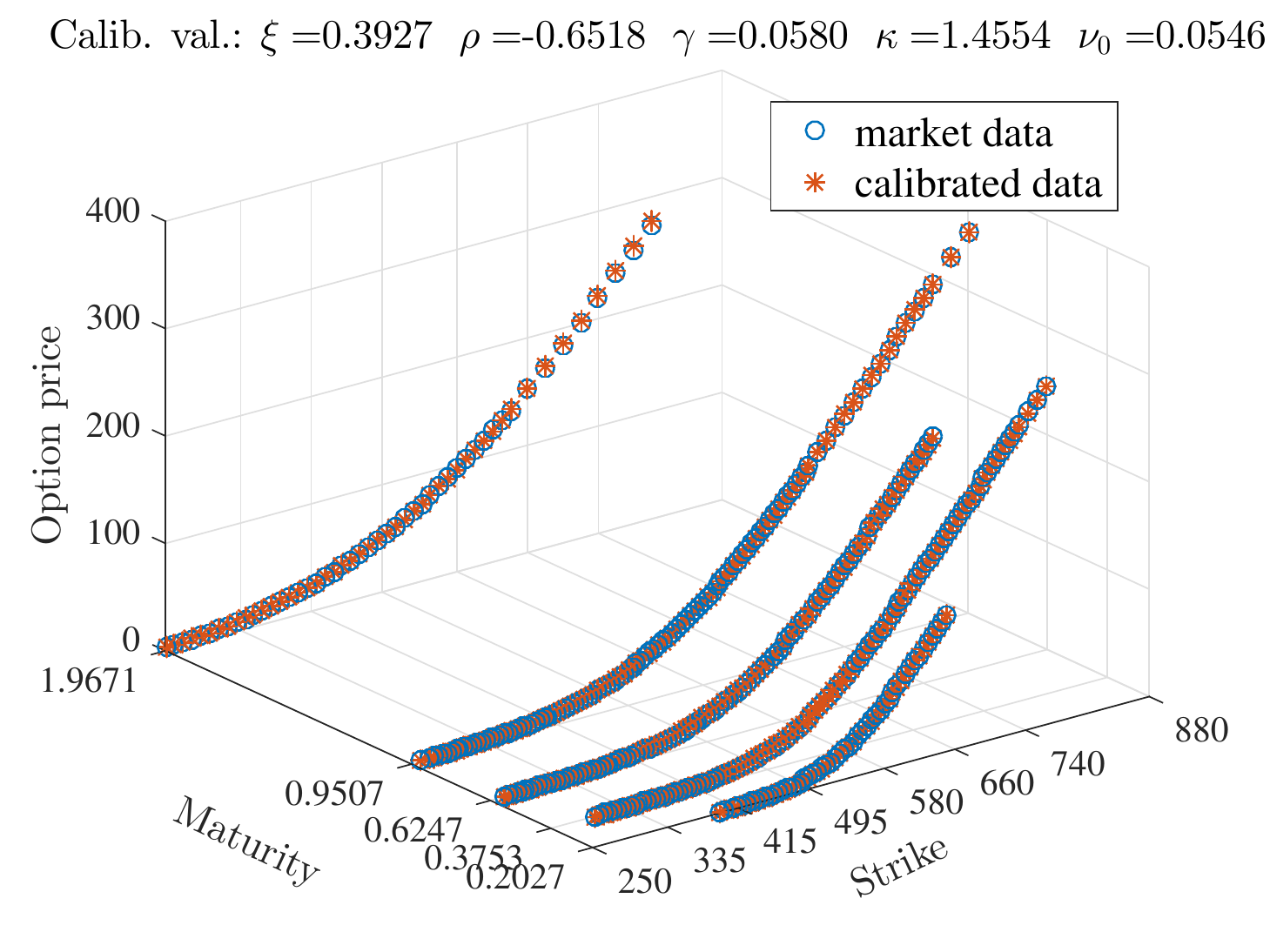}
\includegraphics[width=.45\linewidth]
{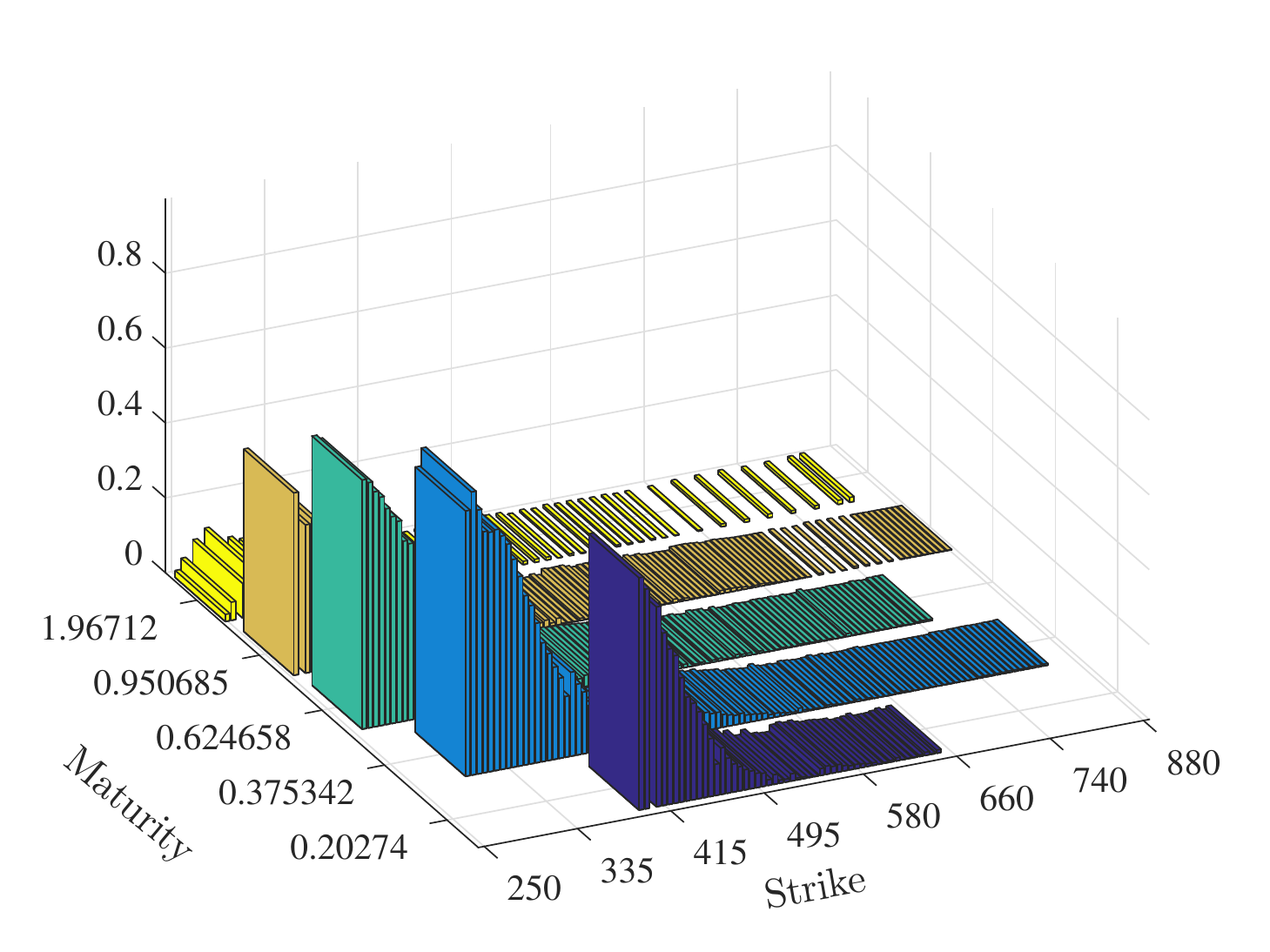}
\centering{De-Americanized problem with closed-form solutions,  $\deam{J}_{\rm CF}(\muopt)$}
 \end{minipage}
 \caption{Left: the Google data set of the de-Americanized American put
options (circles) and the calibrated model data in the Heston model (stars). 
Right: the relative error of
the market and calibrated data, $|P_i^{\rm obs}-P_i^{s,\rm Eu}(\muopt^\star)|/P_i^{\rm obs}$, $i=1,\dots,M$, $s=\{\N,\rm CF\}$} 
\label{ch6:fig_google_ao_deam}
\end{figure}

Figure~\ref{ch6:fig_google_ao_actual} shows the calibration results based on market data for  the detailed and the RBM approaches and the results for the two DAS are provided in  Figure~\ref{ch6:fig_google_ao_deam}.
 The relative error
for all approaches does not exceed $50\%$ and increases in the out-of-the money
region, which corresponds to the smallest option prices.
To reduce this effect,
one could consider different weights in the objective functional, e.g., imposing
larger weights for small option price values.

%

% C O N C L U S I O N

\section{Conclusion}
{In this paper, we applied the reduced basis methodology to the calibration of European and American put options. In the case of American options we additionally considered  a de-Americanization strategy. Both reduction strategies are compared numerically. While RBM techniques aim to achieve smaller dimensions for the discrete spaces in the variational formulations of the problem, the DAS replaces the constrained PDE model of an American option with the unconstrained model of a European option. By doing so, a model error of fixed size occurs, but the RBM offer flexibility by allowing us to adaptively adjust the accuracy.}

% B I B L I O G R A P H Y
\bibliographystyle{siamplain}
\bibliography{references}

\begin{thebibliography}{10}

\bibitem{MR2137495}
{\sc Y.~Achdou}, {\em An inverse problem for a parabolic variational inequality
  arising in volatility calibration with {A}merican options}, SIAM J. Control
  Optim., 43 (2005), pp.~1583--1615 (electronic),
  \href{http://dx.doi.org/10.1137/S0363012903424423}
  {doi:10.1137/S0363012903424423}.

\bibitem{achdou}
{\sc Y.~Achdou and O.~Pironneau}, {\em Computational methods for option
  pricing}, vol.~30 of Frontiers in Applied Mathematics, Society for Industrial
  and Applied Mathematics (SIAM), Philadelphia, PA, 2005,
  \href{http://dx.doi.org/10.1137/1.9780898717495}
  {doi:10.1137/1.9780898717495}.

\bibitem{pir_calib}
{\sc Y.~Achdou and O.~Pironneau}, {\em Numerical procedure for calibration of
  volatility with {A}merican options}, Applied Mathematical Finance, 12 (2005),
  pp.~201--241, \href{http://dx.doi.org/10.1080/1350486042000297252}
  {doi:10.1080/1350486042000297252}.

\bibitem{Amsallem}
{\sc M.~Balajewicz, D.~Amsallem, and C.~Farhat}, {\em Projection-based model
  reduction for contact problems}, 2016,
  \href{http://dx.doi.org/10.1002/nme.5135} {doi:10.1002/nme.5135}.

\bibitem{Balajewicz2016}
{\sc M.~Balajewicz and J.~Toivanen}, {\em Reduced order models for pricing
  {A}merican options under stochastic volatility and jump-diffusion models},
  Procedia Computer Science, 80 (2016), pp.~734--743,
  \href{http://dx.doi.org/10.1016/j.procs.2016.05.360}
  {doi:10.1016/j.procs.2016.05.360}.
\newblock International Conference on Computational Science 2016, {ICCS} 2016,
  6-8 June 2016, San Diego, California, {USA}.

\bibitem{blackscholes}
{\sc F.~Black and M.~Scholes}, {\em The pricing of options and corporate
  liabilities}, J. Polit. Econ., 81 (1973), pp.~637--654,
  \href{http://dx.doi.org/10.1086/260062} {doi:10.1086/260062}.

\bibitem{bouchev}
{\sc I.~Bouchouev and V.~Isakov}, {\em The inverse problem of option pricing},
  Inverse Problems, 13 (1997), pp.~L11--L17,
  \href{http://dx.doi.org/10.1088/0266-5611/13/5/001}
  {doi:10.1088/0266-5611/13/5/001}.

\bibitem{brezzi_raviart}
{\sc F.~Brezzi, W.~W. Hager, and P.-A. Raviart}, {\em Error estimates for the
  finite element solution of variational inequalities. {II}. {M}ixed methods},
  Numer. Math., 31 (1978/79), pp.~1--16,
  \href{http://dx.doi.org/10.1007/BF01396010} {doi:10.1007/BF01396010}.

\bibitem{broadie1997pricing}
{\sc M.~Broadie and P.~Glasserman}, {\em Pricing {A}merican-style securities
  using simulation}, J. Econom. Dynam. Control, 21 (1997), pp.~1323--1352,
  \href{http://dx.doi.org/10.1016/S0165-1889(97)00029-8}
  {doi:10.1016/S0165-1889(97)00029-8}.

\bibitem{burkovskaphd}
{\sc O.~Burkovska}, {\em Reduced Basis Methods for Option Pricing and
  Calibration}, PhD thesis, Technische Universit{\"a}t M{\"u}nchen, 2016.

\bibitem{deam2016}
{\sc O.~Burkovska, K.~Glau, M.~Ga{\ss}, M.~Mahlstedt, W.~Schoutens, and
  B.~Wohlmuth}, {\em Calibration to {A}merican options: Numerical investigation
  of the de-{A}mericanization method}.
\newblock Working paper, in preparation, 2016.

\bibitem{burkovska}
{\sc O.~Burkovska, B.~Haasdonk, J.~Salomon, and B.~Wohlmuth}, {\em Reduced
  basis methods for pricing options with the {B}lack-{S}choles and {H}eston
  models}, SIAM J. Financial Math., 6 (2015), pp.~685--712,
  \href{http://dx.doi.org/10.1137/140981216} {doi:10.1137/140981216}.

\bibitem{carr2010stock}
{\sc P.~Carr and L.~Wu}, {\em Stock options and credit default swaps: A joint
  framework for valuation and estimation}, Journal of Financial Econometrics, 8
  (2010), pp.~409--449, \href{http://dx.doi.org/10.1093/jjfinec/nbp010}
  {doi:10.1093/jjfinec/nbp010}.

\bibitem{vix}
{\sc CBOE}, {\em The {CBOE} volatility index -- {VIX}}, CBOE,  (2009).

\bibitem{clarke}
{\sc N.~Clarke and K.~Parrott}, {\em Multigrid for {A}merican option pricing
  with stochastic volatility}, Applied Mathematical Finance, 6 (1999),
  pp.~177--195, \href{http://dx.doi.org/10.1080/135048699334528}
  {doi:10.1080/135048699334528}.

\bibitem{pironneau}
{\sc R.~Cont, N.~Lantos, and O.~Pironneau}, {\em A reduced basis for option
  pricing}, SIAM J. Financial Math., 2 (2011), pp.~287--316,
  \href{http://dx.doi.org/10.1137/10079851X} {doi:10.1137/10079851X}.

\bibitem{cox1975notes}
{\sc J.~Cox}, {\em Notes on option pricing {I}: Constant elasticity of variance
  diffusions}, Unpublished note, Stanford University, Graduate School of
  Business,  (1975).

\bibitem{crr_tree}
{\sc J.~C. Cox, S.~A. Ross, and M.~Rubinstein}, {\em Option pricing: A
  simplified approach}, Journal of Financial Economics, 7 (1979), pp.~229--263,
  \href{http://dx.doi.org/10.1016/0304-405X(79)90015-1}
  {doi:10.1016/0304-405X(79)90015-1}.

\bibitem{haasdonkopt}
{\sc M.~A. Dihlmann and B.~Haasdonk}, {\em Certified {PDE}-constrained
  parameter optimization using reduced basis surrogate models for evolution
  problems}, Comput. Optim. Appl., 60 (2015), pp.~753--787,
  \href{http://dx.doi.org/10.1007/s10589-014-9697-1}
  {doi:10.1007/s10589-014-9697-1}.

\bibitem{during}
{\sc B.~D{\"u}ring and M.~Fourni{\'e}}, {\em High-order compact finite
  difference scheme for option pricing in stochastic volatility models}, J.
  Comput. Appl. Math., 236 (2012), pp.~4462--4473,
  \href{http://dx.doi.org/10.1016/j.cam.2012.04.017}
  {doi:10.1016/j.cam.2012.04.017}.

\bibitem{egger2005}
{\sc H.~Egger and H.~W. Engl}, {\em Tikhonov regularization applied to the
  inverse problem of option pricing: convergence analysis and rates}, Inverse
  Problems, 21 (2005), pp.~1027--1045,
  \href{http://dx.doi.org/10.1088/0266-5611/21/3/014}
  {doi:10.1088/0266-5611/21/3/014}.

\bibitem{FangOosterlee2011}
{\sc F.~Fang and C.~W. Oosterlee}, {\em {A Fourier-based valuation method for
  Bermudan and barrier options under Heston's Model}}, SIAM Journal of
  Financial Mathematics, 2 (2011), pp.~439--463.

\bibitem{figlewski2000leverage}
{\sc S.~Figlewski and X.~Wang}, {\em Is the 'leverage effect' a leverage
  effect?}, Available at SSRN 256109,  (2000),
  \href{http://dx.doi.org/10.2139/ssrn.256109} {doi:10.2139/ssrn.256109}.

\bibitem{fu2001pricing}
{\sc M.~C. Fu, S.~B. Laprise, D.~B. Madan, Y.~Su, and R.~Wu}, {\em Pricing
  {A}merican options: A comparison of {M}onte {C}arlo simulation approaches},
  Journal of Computational Finance, 4 (2001), pp.~39--88,
  \href{http://dx.doi.org/10.21314/JCF.2001.066} {doi:10.21314/JCF.2001.066}.

\bibitem{MR2891932}
{\sc F.~Gerlich, A.~M. Giese, J.~H. Maruhn, and E.~W. Sachs}, {\em Parameter
  identification in financial market models with a feasible point {SQP}
  algorithm}, Comput. Optim. Appl., 51 (2012), pp.~1137--1161,
  \href{http://dx.doi.org/10.1007/s10589-010-9369-8}
  {doi:10.1007/s10589-010-9369-8}.

\bibitem{UG13}
{\sc S.~Glas and K.~Urban}, {\em On noncoercive variational inequalities}, SIAM
  J. Numer. Anal., 52 (2014), pp.~2250--2271,
  \href{http://dx.doi.org/10.1137/130925438} {doi:10.1137/130925438}.

\bibitem{grepl2005reduced}
{\sc M.~A. Grepl}, {\em Reduced-basis approximation and a posteriori error
  estimation for parabolic partial differential equations}, PhD thesis,
  Massachusetts Institute of Technology, 2005.

\bibitem{Ha13}
{\sc B.~Haasdonk}, {\em Convergence rates of the {POD}-greedy method}, ESAIM
  Math. Model. Numer. Anal., 47 (2013), pp.~859--873,
  \href{http://dx.doi.org/10.1051/m2an/2012045} {doi:10.1051/m2an/2012045}.

\bibitem{haasdonk2011mutrain}
{\sc B.~Haasdonk, M.~Dihlmann, and M.~Ohlberger}, {\em A training set and
  multiple bases generation approach for parameterized model reduction based on
  adaptive grids in parameter space}, Math. Comput. Model. Dyn. Syst., 17
  (2011), pp.~423--442, \href{http://dx.doi.org/10.1080/13873954.2011.547674}
  {doi:10.1080/13873954.2011.547674}.

\bibitem{haasdonk08}
{\sc B.~Haasdonk and M.~Ohlberger}, {\em Reduced basis method for finite volume
  approximations of parametrized linear evolution equations}, M2AN Math. Model.
  Numer. Anal., 42 (2008), pp.~277--302,
  \href{http://dx.doi.org/10.1051/m2an:2008001} {doi:10.1051/m2an:2008001}.

\bibitem{HSW12}
{\sc B.~Haasdonk, J.~Salomon, and B.~Wohlmuth}, {\em A reduced basis method for
  parametrized variational inequalities}, SIAM J. Numer. Anal., 50 (2012),
  pp.~2656--2676, \href{http://dx.doi.org/10.1137/110835372}
  {doi:10.1137/110835372}.

\bibitem{hesthaven2015}
{\sc J.~S. Hesthaven, G.~Rozza, and B.~Stamm}, {\em Certified reduced basis
  methods for parametrized partial differential equations}, SpringerBriefs in
  Mathematics, Springer, Cham; BCAM Basque Center for Applied Mathematics,
  Bilbao, 2016, \href{http://dx.doi.org/10.1007/978-3-319-22470-1}
  {doi:10.1007/978-3-319-22470-1}.
\newblock BCAM SpringerBriefs.

\bibitem{heston}
{\sc S.~L. Heston}, {\em A closed-form solution for options with stochastic
  volatility with applications to bond and currency options}, The Review of
  Financial Studies, 6 (1993), pp.~327--43,
  \href{http://dx.doi.org/10.1093/rfs/6.2.327} {doi:10.1093/rfs/6.2.327}.

\bibitem{hilber}
{\sc N.~Hilber, O.~Reichmann, C.~Schwab, and C.~Winter}, {\em Computational
  methods for quantitative finance. Finite element methods for derivative
  pricing.}, Springer Finance, 2013,
  \href{http://dx.doi.org/10.1007/978-3-642-35401-4}
  {doi:10.1007/978-3-642-35401-4}.

\bibitem{hintermuller}
{\sc M.~Hinterm{\"u}ller}, {\em Inverse coefficient problems for variational
  inequalities: optimality conditions and numerical realization}, M2AN Math.
  Model. Numer. Anal., 35 (2001), pp.~129--152,
  \href{http://dx.doi.org/10.1051/m2an:2001109} {doi:10.1051/m2an:2001109}.

\bibitem{hull03}
{\sc J.~C. Hull}, {\em Options, Futures, and Other Derivative Securities},
  Prentice-Hall, fifth~ed., 2003.

\bibitem{ito2000}
{\sc K.~Ito and K.~Kunisch}, {\em Optimal control of elliptic variational
  inequalities}, Appl. Math. Optim., 41 (2000), pp.~343--364,
  \href{http://dx.doi.org/10.1007/s002459911017} {doi:10.1007/s002459911017}.

\bibitem{janek2011fx}
{\sc A.~Janek, T.~Kluge, R.~Weron, and U.~Wystup}, {\em F{X} smile in the
  {H}eston model}, in Statistical tools for finance and insurance, Springer,
  Heidelberg, 2011, pp.~133--162,
  \href{http://dx.doi.org/10.1007/978-3-642-18062-0\_4}
  {doi:10.1007/978-3-642-18062-0\_4}.

\bibitem{kikuchi}
{\sc N.~Kikuchi and J.~T. Oden}, {\em Contact problems in elasticity: a study
  of variational inequalities and finite element methods}, vol.~8 of SIAM
  Studies in Applied Mathematics, Society for Industrial and Applied
  Mathematics (SIAM), Philadelphia, PA, 1988,
  \href{http://dx.doi.org/10.1137/1.9781611970845}
  {doi:10.1137/1.9781611970845}.

\bibitem{Lee_Seung}
{\sc D.~D. Lee and H.~S. Seung}, {\em {Learning the parts of objects by
  non-negative matrix factorization}}, Nature, 401 (1999), pp.~788--791,
  \href{http://dx.doi.org/10.1038/44565} {doi:10.1038/44565}.

\bibitem{Levendorskii2004}
{\sc S.~Z. Levendorski\u{i}}, {\em Early exercise boundary and option prices in
  {L}\'evy driven models}, Quantitative Finance, 4 (2004), pp.~525--547.

\bibitem{longstaff2001valuing}
{\sc F.~A. Longstaff and E.~S. Schwartz}, {\em Valuing american options by
  simulation: a simple least-squares approach}, Review of Financial studies, 14
  (2001), pp.~113--147, \href{http://dx.doi.org/10.1093/rfs/14.1.113}
  {doi:10.1093/rfs/14.1.113}.

\bibitem{MU14}
{\sc A.~Mayerhofer and K.~Urban}, {\em A reduced basis method for parabolic
  partial differential equations with parameter functions and application to
  option pricing}.
\newblock Preprint, University of Ulm, 2014,
  \href{http://arxiv.org/abs/1408.2709} {arXiv:1408.2709}.
\newblock accepted for publication in Journal of Computational Finance.

\bibitem{mrazek2014optimization}
{\sc M.~Mr{\'a}zek, J.~Posp{\'\i}{\v{s}}il, and T.~Sobotka}, {\em On
  optimization techniques for calibration of stochastic volatility models}, in
  AMCM 2015, November 28--30, 2014, Athens, Greece, 2014, pp.~34--40,
  \url{http://www.inase.org/library/2014/athens/bypaper/APPLIED/APPLIED-04.pdf}.

\bibitem{patera2006}
{\sc A.~T. Patera and G.~Rozza}, {\em Reduced basis approximation and a
  posteriori error estimation for parametrized partial differential equations},
  Tech. Report Version 1.0, MIT 2006--2007, to appear in (tentative rubric) MIT
  Pappalardo Graduate Monographs in Mechanical Engineering, Massachusetts
  Institute of Technology, 2006,
  \url{http://augustine.mit.edu/methodology/bookParts/Patera_Rozza_bookPartI_BV1.pdf}.

\bibitem{red_bs15}
{\sc B.~Peherstorfer, P.~G{\'o}mez, and H.-J. Bungartz}, {\em Reduced models
  for sparse grid discretizations of the multi-asset {B}lack-{S}choles
  equation}, Advances in Computational Mathematics, 41 (2015), pp.~1365--1389,
  \href{http://dx.doi.org/10.1007/s10444-015-9421-4}
  {doi:10.1007/s10444-015-9421-4}.

\bibitem{pironneau09}
{\sc O.~Pironneau}, {\em Calibration of options on a reduced basis}, J. Comput.
  Appl. Math., 232 (2009), pp.~139--147,
  \href{http://dx.doi.org/10.1016/j.cam.2008.10.070}
  {doi:10.1016/j.cam.2008.10.070}.

\bibitem{pironneau11}
{\sc O.~Pironneau}, {\em Reduced basis for vanilla and basket options}, Risk
  and Decision Analysis, 2 (2011), pp.~185--194,
  \href{http://dx.doi.org/10.3233/RDA-2011-0045} {doi:10.3233/RDA-2011-0045}.

\bibitem{pironneau12}
{\sc O.~Pironneau}, {\em Proper orthogonal decomposition for pricing options},
  Journal of Computational Finance, 16 (2012), pp.~33--46,
  \href{http://dx.doi.org/10.21314/JCF.2012.246} {doi:10.21314/JCF.2012.246}.

\bibitem{quarteroni_rbm}
{\sc A.~Quarteroni, A.~Manzoni, and F.~Negri}, {\em Reduced basis methods for
  partial differential equations}, vol.~92 of Unitext, Springer, Cham, 2016,
  \href{http://dx.doi.org/10.1007/978-3-319-15431-2}
  {doi:10.1007/978-3-319-15431-2}.
\newblock An introduction, La Matematica per il 3+2.

\bibitem{QV}
{\sc A.~Quarteroni and A.~Valli}, {\em Numerical approximation of partial
  differential equations}, vol.~23 of Springer Series in Computational
  Mathematics, Springer-Verlag, Berlin, 1994.

\bibitem{rogers2002monte}
{\sc L.~C.~G. Rogers}, {\em Monte {C}arlo valuation of {A}merican options},
  Math. Finance, 12 (2002), pp.~271--286,
  \href{http://dx.doi.org/10.1111/1467-9965.02010}
  {doi:10.1111/1467-9965.02010}.

\bibitem{rovas}
{\sc D.~Rovas}, {\em Reduced-basis output bound methods for parametrized
  partial differential equations}, PhD thesis, Massachusetts Institute of
  Technology, 2003.

\bibitem{rozza_2007}
{\sc G.~Rozza and K.~Veroy}, {\em On the stability of the reduced basis method
  for {S}tokes equations in parametrized domains}, Comput. Methods Appl. Mech.
  Engrg., 196 (2007), pp.~1244--1260,
  \href{http://dx.doi.org/10.1016/j.cma.2006.09.005}
  {doi:10.1016/j.cma.2006.09.005}.

\bibitem{rubinstein1994implied}
{\sc M.~Rubinstein}, {\em Implied binomial trees}, The Journal of Finance, 49
  (1994), pp.~771--818, \href{http://dx.doi.org/10.2307/2329207}
  {doi:10.2307/2329207}.

\bibitem{sachs2014}
{\sc E.~W. Sachs and M.~Schneider}, {\em Reduced-order models for the implied
  variance under local volatility}, Int. J. Theor. Appl. Finance, 17 (2014),
  pp.~1450053, 23, \href{http://dx.doi.org/10.1142/S0219024914500538}
  {doi:10.1142/S0219024914500538}.

\bibitem{sachs2014_1}
{\sc E.~W. Sachs, M.~Schneider, and M.~Schu}, {\em Adaptive trust-region {POD}
  methods in {PIDE}-constrained optimization}, in Trends in {PDE} constrained
  optimization, vol.~165 of Internat. Ser. Numer. Math., Birkh\"auser/Springer,
  Cham, 2014, pp.~327--342,
  \href{http://dx.doi.org/10.1007/978-3-319-05083-6\_20}
  {doi:10.1007/978-3-319-05083-6\_20}.

\bibitem{sachs2008}
{\sc E.~W. Sachs and M.~Schu}, {\em Reduced order models ({POD}) for
  calibration problems in finance}, in Numerical Mathematics and Advanced
  Applications: Proceedings of {ENUMATH} 2007, the 7th {E}uropean Conference on
  Numerical Mathematics and Advanced Applications, Graz, Austria, September
  2007, K.~Kunisch, G.~Of, and O.~Steinbach, eds., Springer Berlin Heidelberg,
  Berlin, Heidelberg, 2008, pp.~735--742,
  \href{http://dx.doi.org/10.1007/978-3-540-69777-0\_88}
  {doi:10.1007/978-3-540-69777-0\_88}.

\bibitem{volkwein2010}
{\sc E.~W. Sachs and S.~Volkwein}, {\em P{OD}-{G}alerkin approximations in
  {PDE}-constrained optimization}, GAMM-Mitt., 33 (2010), pp.~194--208,
  \href{http://dx.doi.org/10.1002/gamm.201010015} {doi:10.1002/gamm.201010015}.

\bibitem{schiela}
{\sc A.~Schiela and D.~Wachsmuth}, {\em Convergence analysis of smoothing
  methods for optimal control of stationary variational inequalities with
  control constraints}, ESAIM Math. Model. Numer. Anal., 47 (2013),
  pp.~771--787, \href{http://dx.doi.org/10.1051/m2an/2012049}
  {doi:10.1051/m2an/2012049}.

\bibitem{Schoutens2004}
{\sc W.~Schoutens, E.~Simons, and J.~Tistaert}, {\em {A perfect calibration!
  Now what?}}, Wilmott Magazine,  (2004), pp.~66--78.

\bibitem{vexler}
{\sc B.~Vexler}, {\em Adaptive Finite Element Methods for Parameter
  Identification Problems}, PhD thesis, University of Heidelberg, 2004.

\bibitem{Woh00a}
{\sc B.~I. Wohlmuth}, {\em A mortar finite element method using dual spaces for
  the {L}agrange multiplier}, SIAM J. Numer. Anal., 38 (2000), pp.~989--1012,
  \href{http://dx.doi.org/10.1137/S0036142999350929}
  {doi:10.1137/S0036142999350929}.

\bibitem{VeroyParisWorkshop}
{\sc Z.~Zhang, E.~Bader, and K.~Veroy}, {\em A slack approach to reduced-basis
  approximation and error estimation for variational inequalities}, C. R. Math.
  Acad. Sci. Paris, 354 (2016), pp.~283--289,
  \href{http://dx.doi.org/10.1016/j.crma.2015.10.024}
  {doi:10.1016/j.crma.2015.10.024}.

\end{thebibliography}

\newpage
\appendix
\section{Google market data}
\begin{table}[h]
\begin{tiny}\begin{tabular}{|c|c|c|c|c|c|}
\hline
$K\backslash T$&0.2027 &  0.3753&   0.6247&  0.9507& 1.9671\\
\hline\hline
250& & & & &2.50\\
260& & & &1.20&2.90\\
265& & &0.55& & \\
270& & &0.62&1.30&3.30\\
275& & &0.68&1.50& \\
280& & &0.75&1.62&3.80\\
285& & &0.80&1.73& \\
290& & &0.88&1.80&4.65\\
295& & &0.97&1.95& \\
300& &0.28&1.00&2.10&5.35\\
305& &0.40&1.12&2.27& \\
310& &0.40&1.20&2.45&6.25\\
315& &0.40&1.30&2.67& \\
320& &0.47&1.38&2.85&7.40\\
325& &0.57&1.50&3.12& \\
330& &0.65&1.60&3.40&8.10\\
335& &0.72&1.75&3.60& \\
340& &0.78&1.90&3.85&9.65\\
345& &0.82&2.05&4.20& \\
350& &0.88&2.23&4.55&10.90\\
355& &0.97&2.40&4.90& \\
360& &1.05&2.70&5.30&12.55\\
365& &1.12&2.90&5.65& \\
370& &1.25&3.17&6.15&14.30\\
375& &1.38&3.45&6.65& \\
380& &1.52&3.75&7.15&16.45\\
385& &1.65&4.10&7.65& \\
390& &2.05&4.45&8.15&18.30\\
395&0.93&2.25&4.85&8.80& \\
400&1.05&2.50&5.30&9.45&20.85\\
405& &2.73&5.75&10.10& \\
410&1.40&3.08&6.30&10.75&22.85\\
415&1.48&3.35&6.85&11.55& \\
420&1.65&3.75&7.45&12.40&25.55\\
425&1.93&4.10&8.15&13.40& \\
430&2.10&4.55&8.80&14.20&28.30\\
435&2.33&5.05&9.55&15.35& \\
440&2.70&5.55&10.40&16.25&31.35\\
445&3.10&6.10&11.25&17.35& \\
450&3.50&6.80&12.30&18.65&34.00\\
455&3.92&7.50&13.25&19.80& \\
460&4.25&8.35&14.25&21.20&37.80\\
465&5.20&9.20&15.45&22.60& \\
470&5.85&10.20&16.70&23.90&41.45\\
475&6.65&11.25&18.00&25.30& \\
480&7.60&12.40&19.45&27.15&45.20\\
485&8.65&13.70&20.95&28.55& \\
490&9.85&15.10&22.55&30.75&49.20\\
495&11.15&16.55&24.25&32.60& \\
500&12.65&18.25&26.00&34.60&53.50\\
505&14.05&20.00&27.90&36.80& \\
510&16.05&21.95&30.00&38.70&58.15\\
515&17.75&23.90&32.15&41.45& \\
520&19.90&26.20&34.40&43.65&62.90\\
525&22.50&28.55&36.70&46.05& \\
530&24.70&31.05&39.20&48.55&67.90\\
535&27.40&33.50&41.85&50.85& \\\hline
\end{tabular}
\end{tiny}
\begin{tiny}\begin{tabular}{|c|c|c|c|c|c|}
\hline
$K\backslash T$&0.2027 &  0.3753&   0.6247&  0.9507& 1.9671\\
\hline\hline
540&30.20&36.40&44.45&52.70&72.65\\
545&33.25&39.25&47.60&56.10& \\
550&37.15&42.30&50.25&58.95&78.50\\
555&40.80&45.45&53.15&61.85& \\
560&44.30&49.00&56.55&64.80&84.05\\
565&48.25&52.35&59.40&67.50& \\
570&51.95&56.00&63.50&71.10&90.40\\
575&55.90&59.80&66.70&73.90& \\
580&60.20&63.70&69.85&77.25&96.45\\
585&64.50&67.10&74.25&81.30& \\
590&68.80&71.60&77.85&84.45&103.00\\
595&73.35&76.30&81.70&88.20& \\
600&77.35&80.65&85.10&91.55&108.85\\
605&82.35&84.80&89.65&95.20& \\
610&87.30&89.65&93.20&99.10&116.45\\
615&92.50&94.05&97.40&102.90& \\
620&96.65&98.45&101.85&106.80&123.75\\
625&101.30&103.00&106.40&111.80& \\
630&106.50&107.25&110.60&116.00&130.65\\
635&111.75&112.95&115.00&120.25& \\
640&117.20&117.25&119.85&124.40&138.10\\
645&121.45&122.30&124.40&128.55& \\
650&127.25&126.70&128.90&132.55& \\
655& &131.45&133.70&137.15& \\
660& &136.45&138.15&141.80&153.70\\
665& &141.45&143.05&146.10& \\
670& &146.65&147.35&150.50& \\
675& &151.80&152.05&154.90& \\
680& &156.20&158.25&159.45&169.45\\
685& &161.20&161.60&164.05& \\
690& &166.30&166.55&168.65& \\
695& &171.30&171.45&173.15& \\
700& &176.20&176.45&178.15&186.40\\
705& &181.10&181.75& & \\
710& &186.15&186.25&188.40& \\
715& &191.10&191.35& & \\
720& &196.05&195.85&197.15&203.70\\
725& &201.05&201.80& & \\
730& &205.85&206.50&206.95& \\
735& &211.00&211.05& & \\
740& &216.80&216.80&216.20&222.30\\
745& &221.80&221.90& & \\
750& &226.55&226.85&226.20& \\
755& &231.80& & & \\
760& &236.80& &236.15&240.60\\
765& &241.95& & & \\
770& &246.95& &246.05& \\
775& &251.80& & & \\
780& &256.80& &256.05&258.90\\
790& &267.00& &265.70& \\
800& &276.95& &276.05&277.70\\
810& &287.25& &286.05& \\
820& & & &296.05& \\
830& & & &306.05& \\
840& & & &316.05& \\
860& & & &336.05& \\
880& & & &356.05& \\\hline
\end{tabular}
\end{tiny}
\caption{Google market data consisting of 401 American put options with $S_0=523,755$ on February 2nd, 2015.}\label{ch6:table:google_data}
\end{table}

\end{document}